%% file: Royset_monit_rev.tex
\documentclass[titlepage,final,11pt]{article}
\input mac.tex

\begin{document}


\begin{center}
\begin{large}
{\bf Variational Analysis of a Nonconvex and Nonsmooth Optimization Problem:\\
An Introduction}
\smallskip
\end{large}
\vglue 0.5truecm
\begin{tabular}{cc}
  \begin{large} {\sl Johannes O. Royset} \end{large}\\
  Daniel J. Epstein Department of Industrial and Systems Engineering\\
  University of Southern California
\end{tabular}

\vskip 1.5truecm

{\em Dedicated to Boris Mordukhovich for his fundamental and lasting contributions to variational analysis}

\end{center}

\vskip 1.5truecm

\noindent {\bf Abstract}. Variational analysis provides the theoretical foundations and practical tools for constructing optimization algorithms without being restricted to smooth or convex problems. We survey the central concepts in the context of a concrete but broadly applicable problem class from composite optimization in finite dimensions. While prioritizing accessibility over mathematical details, we introduce subgradients of arbitrary functions and the resulting optimality conditions, describe approximations and the need for going beyond pointwise and uniform convergence, and summarize proximal methods. We derive dual problems from parametrization of the actual problem and the resulting relaxations. The paper ends with an introduction to second-order theory and its role in stability analysis of optimization problems. 

\vskip 0.4truecm

\halign{&\vtop{\parindent=0pt
   \hangindent2.5em\strut#\strut}\cr
                         \cr

{\bf Date}:\quad \ \today \cr}

\baselineskip=15pt

\section{Introduction}\label{sec:intro}

Finding a minimizer of a nonconvex function is fundamentally hard. The challenge is exacerbated when the function is nonsmooth because the traditional approach of computing a point with vanishing gradient is no longer available. In parallel, the shift from problems with equality constraints defined by smooth functions to those with more general constraints implicitly introduces nonsmoothness that invalidates the classical approach of Lagrange. Variational analysis grew from convex analysis and calculus of variation to address nonconvex and nonsmooth problems. It provides the fundamental mathematical tools for analyzing such problems, characterizing the solutions, and justifying their algorithms. Variational analysis is undeterred by nonsmoothness and nonconvexity, but can leverage these properties when present in some part of a problem. 

We survey the main tools of variational analysis in the context of a specific class of optimization problems with many applications. For $X\subset\reals^n$ and $G:\reals^n\to \reals^m$, we consider the problem  
\begin{equation}\label{eqn:actualproblem}
\nnmin_{x\in X} ~h\big(G(x)\big),
\end{equation}
where $h:\reals^m\to \Reals = [-\infty, \infty]$ is an {\em extended real-valued} function. The fact that functions can take the values $-\infty$ and $\infty$ is a feature that distinguishes variational analysis from most other areas of mathematics. It provides a convenient way of encoding constraints using infinite penalties and thus unifies the treatment of constrained optimization problems with those lacking such restrictions. Specifically, we can reformulate \eqref{eqn:actualproblem} as the problem of minimizing a function $\phi$ defined as $\phi(x) = h(G(x))$ when $x\in X$ and $\phi(x) = \infty$ otherwise. 

Regardless of the formulation, we bring forth the three components $X$, $h$, and $G$ to stress the importance of problem structure in nonconvex and nonsmooth optimization. If nothing is known about a function to be minimized, then there is little room for variational analysis and it becomes difficult to beat the naive solution strategy of randomly trying different points. However, in the presence of problem structure powerful techniques and rules from variational analysis emerge. Throughout the paper, we assume that 
\[
\mbox{$X$ is a {\em nonempty}, {\em closed}, and {\em convex} set and $G$ is a {\em smooth} (i.e., $C^1$) mapping.} 
\]
While many of the facts stated below hold more generally, there are specific advantages to these assumptions that allow us to develop simpler expressions. Moreover, for some {\em nonempty polyhedral} set $Y\subset \reals^m$ and a {\em symmetric positive semidefinite} $m\times m$-matrix $Q$, we universally assume that 
\begin{equation}\label{eqn:h}
h(z) = \sup_{y\in Y} ~\langle y,z\rangle - \tfrac{1}{2}\langle y, Qy\rangle~~~~~\forall z\in\reals^m.
\end{equation}
At first, the specific form of $h$ may appear limiting but many examples in \cite{AravkinBurkePillonetto.14,DavisDrusvyatskiyPaquette.20,CharisopoulosDavisDiazDrusvyatskiy.21,FerrisHuberRoyset.24} and below illustrate its versatility. Applications include nonlinear optimization, phase retrieval, robust estimation, Kalman smoothing, sensing, risk management, and most significantly, stochastic optimization where the function class appeared originally \cite{RockafellarWets.86}. The function $h$ is convex and in fact {\em piecewise linear-quadratic} \cite[Section 10.E, Example 11.18]{VaAn}; see also \cite{Rockafellar.87,Sun.92}. While an introductory course on optimization may refer to convex functions as ``easy,'' the situation is more nuanced and $h$ represents a class of {\em computationally attractive} convex functions. Under the stated assumptions, \eqref{eqn:actualproblem} may lack both convexity and smoothness overall as seen from the example $X = \reals$, $G(x) = x^2-1$, $h(z) = |z| = \sup_{y\in [-1,1]} yz$. Nevertheless, the components have such properties.

After illustrating the breadth of the problem class \eqref{eqn:actualproblem} in Section \ref{sec:examples} using five examples, we start the survey with basic concepts in Section \ref{sec:basic}. These include the extension from gradients to subgradients, and from normal subspaces to normal cones. Section \ref{sec:firstorderoptimality} puts the extensions to use for the purpose of generalizing the classical necessary condition for a minimizer of a smooth function: vanishing gradient. Section \ref{sec:approx} examines the effect of approximations and shows that the usual notions of pointwise and uniform convergence need to be supplemented when examining optimization problems. Section \ref{sec:algo} briefly discusses algorithms for \eqref{eqn:actualproblem}. Section \ref{sec:duality} introduces parametrization and the resulting relaxations and dual problems. The paper ends in Section \ref{sec:secondorder} with a glimpse at second-order theory, which extends Hessian matrices to the nonsmooth setting.  

Throughout, we focus on concepts and terminology, but give many references to sources containing proofs and further reading. Additional examples and motivation appear in the tutorials \cite{Royset.21,FerrisHuberRoyset.24}. We only briefly touch on the history of variational analysis and refer to the commentaries in \cite{VaAn,Mordukhovich.18}. Essentially all the concepts discussed in this introductory survey have generalizations beyond finite dimensions; see, for example, \cite{Mordukhovich.13a,Mordukhovich.13b}.

\section{Examples}\label{sec:examples}

The problem class \eqref{eqn:actualproblem} arises naturally in applications as seen from five examples. 

\begin{example}{\rm (goal optimization).}\label{ex:goal} In multi-objective optimization, we might hope to minimize the smooth functions $g_i:\reals^n\to \reals$, $i=1, \dots, m$. This is rarely possible, however, because no single $x$ minimizes all the functions. An approach to explore the trade-off between these functions is goal optimization, where one aims to determine a solution $x \in X$ such that each $g_i(x)$ is not higher than a target value $\tau_i$. If the target value is not reached, the solution receives a per-unit penalty of $\alpha_i \in [0,\infty)$. For nonempty, closed, and convex $X\subset\reals^n$, this leads to the formulation
\[
\nnmin_{x\in X} \sum_{i=1}^m \alpha_i \max\big\{0, g_i(x) - \tau_i\big\},
\]
which is of the form \eqref{eqn:actualproblem} with $G(x) = (g_1(x)-\tau_1, \dots, g_m(x)-\tau_m)$ and 
\[
h(z) = \sup_{y\in Y} \, \langle y,z\rangle, ~\mbox{ where } Y = [0,\alpha_1] \times \cdots \times [0,\alpha_m].
\]
Here, $h$ is real-valued because $Y$ is bounded.
\end{example}

\begin{example}{\rm (nonlinear optimization).}\label{ex:nonlinear} The smooth functions $g_i:\reals^n\to \reals$, $i = 0, 1, \dots, m+q$ produce the nonlinear optimization problem 
\[
\nnmin_{x\in\reals^n} \, g_0(x)~ \mbox{ subject to } ~g_i(x) = 0, ~i=1, \dots, m; ~~g_{i}(x) \leq 0, ~i=m+1, \dots, m+q.  
\]
It fits the form \eqref{eqn:actualproblem} with $X = \reals^n$, $G(x) = (g_0(x), g_1(x), \dots, g_{m+q}(x))$, $Y = \{1\} \times \reals^m \times [0, \infty)^q $, and for $z = (z_0, z_1, \dots, z_{m+q})$, 
\[
h(z) = \sup_{y\in Y} \, \langle y, z\rangle = \begin{cases}
  z_0 & \mbox{ if } z_i = 0 ~\mbox{ for } i =1, \dots, m;~ z_i \leq 0 \mbox{ for } i = m+1, \dots, m+q\\
  \infty & \mbox{ otherwise}.   
\end{cases}
\]
Thus, $h$ assigns an infinite penalty to any $x$ that violates the stated equality and inequality constraints. 
\end{example}

\begin{example}{\rm (risk minimization).}\label{ex:risk} We may seek a decision $x\in X$ that is good across $m$ different future scenarios. If the cost of a decision in the $i$th scenario is described by the smooth function $g_i:\reals^n\to \reals$ and $p_i$ is the probability of that scenario, then one could consider the average cost: $\sum_{i=1}^m p_i g_i(x)$. However, a risk-averse decision maker might prefer to consider a worst-case cost in some sense. For $\alpha \in [0,1]$, this could result in a formulation that minimizes the $\alpha$-superquantile (a.k.a. CVaR, AVaR, expected shortfall) across the scenarios; see, e.g., \cite{Royset.24}. Specifically, one obtains the formulation 
\[
\nnmin_{x\in X} \, \sup_{y\in Y} \sum_{i=1}^m y_i g_i(x),~~\mbox{where } Y = \Big\{y\in [0,\infty)^m~\Big|~\sum_{i=1}^m y_i = 1,  ~(1-\alpha)y_i \leq p_i~\forall i\Big\}.
\]
If $X$ is a nonempty, closed, and convex set, then the formulation fits the mold of \eqref{eqn:actualproblem} with $G(x) = (g_1(x), \dots, g_m(x))$, and $h(z) = \sup_{y\in Y} \,\langle y, z\rangle$. In this case, $h$ is real-valued because $Y$ is bounded.
\end{example}

\begin{example}{\rm (regression and inverse problems).}\label{ex:stat}  Lasso regression seeks to find a sparse affine function on $\reals^n$ that fits a given data set well in a least-square sense. With the data represented by the $m\times n$-matrix $A$ and the vector $b$, this leads to the formulation  
\[
\nnmin_{x\in \reals^n} \, \|Ax - b\|_2^2 + \theta \|x\|_1,
\]
where $\theta \in (0,\infty)$ is a regularization parameter. The formulation fits the form \eqref{eqn:actualproblem} with $G(x) = (g_0(x), g_1(x), \dots, g_n(x))$, $g_0(x) = \|Ax -b\|_2^2$, $g_j(x) = \theta x_j$, $j=1, \dots, n$, and $h(z) = \sup_{y\in Y} \,\langle y, z\rangle$, with $Y = \{1\} \times [-1, 1]^n$. Thus, $h$ is real-valued because $Y$ is bounded. Sometimes the $\|\cdot\|_1$-term is replaced by a nonconvex regularization that tapers off away from zero. One possibility is to replace $g_j$, $j\in \{1, \dots, n\}$, by 
\[
g_{j}(x) = \begin{cases}
  2\theta - \theta \exp(1-x_{j}) & \mbox{ if } x_{j}\in (1,\infty)\\
  \theta x_{j}           & \mbox{ if } x_{j}\in [-1, 1]\\
  \theta \exp(1+x_{j})-2\theta  & \mbox{ otherwise,}\\
\end{cases}
\]
which is smooth. Thus, the problem remains in the form \eqref{eqn:actualproblem}.  

Phase retrieval is an example of an inverse problem where one seeks to reconstruct a signal from recorded Fourier magnitudes. Given squared Fourier magnitudes $b_1, \dots, b_m$, the goal becomes to determine $x$ such that $b_i = |\langle a^i, x\rangle|^2$, where $a^1, \dots, a^m$ are known; see e.g., \cite{DuchiRuan.19}. (For simplicity, we ignore that these vectors as well as $x$ are generally complex-valued.) This leads to the formulation
\[
\nnmin_{x\in\reals^n} \frac{1}{m} \nsum_{i=1}^m  \big| \langle a^i, x\rangle^2 -b_i\big|,
\]
which again is of the form \eqref{eqn:actualproblem}: $G(x)= (\langle a^1, x\rangle^2 -b_1, \dots, \langle a^m, x\rangle^2 -b_m)/m$ and $h(z) = \nsup_{y\in [-1,1]^m} \langle y, z\rangle$.
\end{example}

The form \eqref{eqn:actualproblem} reaches beyond optimization to more general problems such as variational inequalities and their wide range of applications; see, e.g., \cite{FerrisPang.97,FacchineiPang.03} and \cite[Chapter 7]{primer}.

\begin{example}{\rm (variational inequality).}\label{ex:VI} For a smooth mapping $F:\reals^n\to \reals^n$ and a nonempty, closed, and convex set $C\subset \reals^n$, the variational inequality of finding
\begin{equation}\label{eqn:va}
x\in C~~~\mbox{ such that } ~~~\big\langle F(x), y - x \big\rangle \geq 0 ~~~\forall y \in C
\end{equation}
can be solved via \eqref{eqn:actualproblem} because of the following fact (see, e.g., \cite[Example 7.13]{primer}): A point $x^\star$ solves the variational inequality if and only if $x^\star$ minimizes $g$ over $C$ and $g(x^\star)= 0$, where  
\[
g(x) = \sup_{y\in C} \big\langle F(x), x-y\big\rangle = \sup_{(\alpha,y) \in \{1\} \times C} \Big\langle (\alpha,y), \big(\langle F(x),x\rangle, -F(x) \big) \Big\rangle. 
\]
Thus, one can take $X = C$, $Y = \{1\}\times C$, and $G(x) = (\langle F(x),x\rangle, -F(x) )$ in \eqref{eqn:actualproblem}.

Variational inequalities arise, for instance, in the context of equilibrium problems. Example 7.3 of \cite{primer} lays out one possibility stemming from \cite{Samuelson.52,TakayamaJudge.71}. A product is made by $m$ geographically dispersed manufacturer and bought by consumers in $n$ different regions. A variational inequality describes how many units $s_i$ manufacturer $i$ will supply, the number of units $d_j$ that region $j$ will demand, and the quantity $w_{ij}$ that will be transported from $i$ to $j$. With $s = (s_1, \dots, s_m)$ and $d = (d_1, \dots, d_n)$, we assume that $p_i(s)$ is the price by manufacturer $i$ given supply vector $s$ and $q_j(d)$ is the price in region $j$ under demand vector $d$, with these being known functions. With $w = (w_{ij}, i=1, \dots, m, j=1, \dots, n)$ and $c_{ij}(w)$ being the cost of transportation between $i$ and $j$ given $w$, an economic equilibrium condition states that $x=(s,d,w)$ must satisfy the variational inequality \eqref{eqn:va} defined by 
\[
C = \Big\{(s,d,w)\in\reals^m\times\reals^n\times \reals^{mn} ~\Big|~ w_{ij} \geq 0, ~\nsum_{j=1}^n w_{ij} = s_i, ~\nsum_{i=1}^m w_{ij} = d_j ~~ ~\forall i,j\Big\}
\]
and $F(x) = (p_1(s), \dots, p_m(s), -q_1(d), \dots, -q_n(d), c_{11}(w), \dots, c_{mn}(w))$, with $x = (s,d,w)$.
\end{example}

\section{Variational Analysis: Basic Concepts}\label{sec:basic}

We start by introducing the fundamental concepts. For a function $f:\reals^n\to \Reals$, $\dom f = \{x \in \reals^n\,|\,f(x)<\infty\}$ is its {\em domain} and $\epi f = \{(x,\alpha) \in \reals^n \times \reals~|~f(x) \leq \alpha \}$ is its {\em epigraph}. While ``domain'' might in other areas pertain to the allowable arguments for a function, we stress that $f$ is defined on the whole of $\reals^n$ and $\dom f$ simply indicates where the function values are not infinity. The function is {\em lower semicontinuous} (lsc) if $\epi f$ is a closed subset of $\reals^n\times \reals$. It is {\em convex} if $\epi f$ is a convex set and it is {\em proper} if $\epi f \neq \emptyset$ and $f(x) > -\infty$ for all $x\in\reals^n$. In particular, $h$ in \eqref{eqn:h} is proper, lsc, and convex. 

Moreover, $\inf f = \inf_{x\in \reals^n} f(x)$ is the {\em minimum value} of $f$, and $\epsilon$-$\nargmin f = \{x\in \dom f~|~f(x) \leq \inf f + \epsilon\}$ is its set of $\epsilon$-{\em minimizers}, with $\epsilon \in [0,\infty)$. If $\epsilon = 0$, then we simply write $\nargmin f$ and refer to its points as {\em minimizers}. In particular, the function having $f(x) = \infty$ for all $x\in \reals^n$ has no $\epsilon$-minimizers because $\dom f = \emptyset$.  A point $x^\star\in \dom f$ is a {\em local minimizer} of $f$ if there exists $\rho >0$ such that $f(x^\star) \leq f(x)$ for all $x$ satisfying $\|x-x^\star\|_2 \leq \rho$. 

The {\em indicator function} of a set $C$ has $\iota_C(x) = 0$ if $x\in C$ and $\iota_C(x) = \infty$ otherwise. Thus, \eqref{eqn:actualproblem} can simply be written in terms of an extended real-valued function $\phi:\reals^n\to \Reals$:  
\begin{equation}\label{eqn:actualproblem2}
\nnmin_{x\in\reals^n} \, \phi(x) = \iota_X(x) + h\big(G(x)\big). 
\end{equation}
We refer to this problem as the {\em actual problem} and retain the definition of $\phi$ throughout. 

The trivial reformulation from \eqref{eqn:actualproblem} to \eqref{eqn:actualproblem2} is especially useful as we develop extensions of the classical optimality condition of vanishing gradients at local minimizers. Optimality conditions for \eqref{eqn:actualproblem} arise by considering gradients of $\phi$. However, any such development requires extension of gradients to the nonsmoothness case and this leads to subgradients. 

\begin{figure}
\centering
\includegraphics[width=0.80\textwidth]{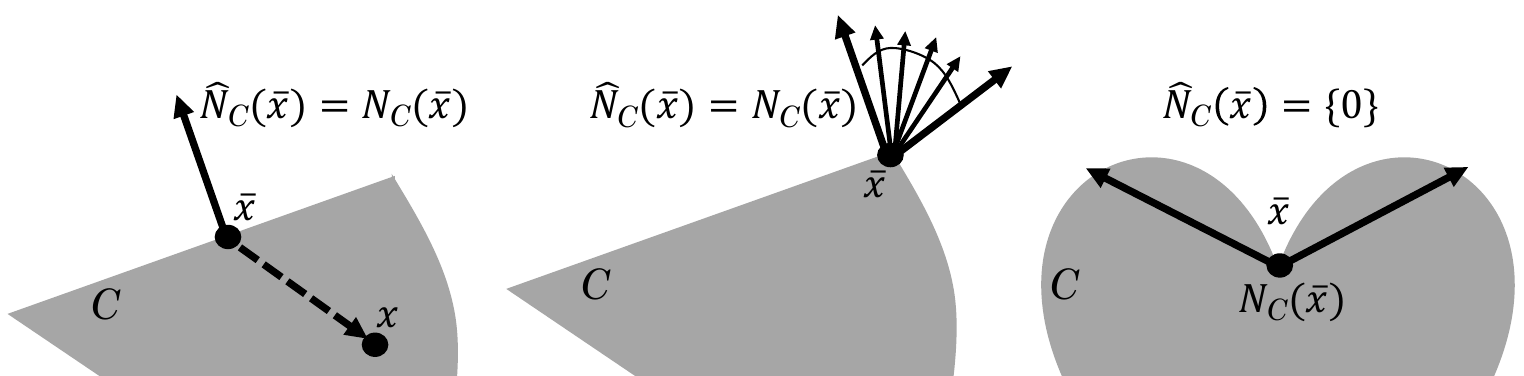}
\caption{Regular normal vectors and normal vectors (solid arrows) at a point of smoothness (left), at a kink (middle), and at a cusp (right).}\label{fig:normals}
\end{figure}

Subgradients of an arbitrary function $f:\reals^n\to \Reals$ are intimately connected to normal vectors and, following the path taken by the pioneering work in \cite{Mordukhovich.76}, we use the latter to defined the former. A normal vector to a smooth manifold at a point on the manifold is well understood and defines a normal subspace. Variational analysis extends this concept to arbitrary sets which may have kinks, cusps, and other irregularities. For a set $C\subset\reals^n$ and a point $\bar x\in C$, a vector $v$ is a {\em regular normal vector} to $C$ at $\bar x$ if 
\[
\limsup_{\substack{x\in C\to \bar x \\ x\neq \bar x}} \frac{\langle v, x-\bar x \rangle}{\|x-\bar x\|_2} \leq 0.
\]
The set of all regular normal vectors to $C$ at $\bar x$ is denoted by $\widehat{N}_C(\bar x)$. Since $\langle v, x-\bar x\rangle$ is nonpositive when the angle between the vectors $v$ and $x-\bar x$ is in the range 90 to 270 degrees, the definition indeed reflects the expected relation between regular normal vectors and vectors ``pointing inwards'' from $\bar x$ see Figure \ref{fig:normals}(left).  

A vector $v$ is {\em normal} to $C$ at $\bar x$ if there are sequences $x^\nu \in C\to \bar x$ and $v^\nu \in \widehat N_C(x^\nu) \to v$. The set of all normal vectors to $C$ at $\bar x$ is the {\em normal cone} to $C$ at $\bar x$ and is  denoted by $N_C(\bar x)$. By convention, we define $\widehat N_C(\bar x) = N_C(\bar x) = \emptyset$ if $\bar x\not\in C$. A normal cone in this sense is sometimes referred to as the Mordukhovich, general, basic, or limiting normal cone. Figure \ref{fig:normals} illustrates regular normal vectors and normal vectors for sets in $\reals^2$ and how there might be more of the latter than the former; see the right portion of the figure where the limits of regular normal vectors at points near $\bar x$ contribute to $N_C(\bar x)$. 

There is an extensive calculus for computing normal cones in various settings; see, e.g., \cite[Theorem 6.14]{VaAn} and \cite[Theorem 4.46]{primer}. We recall a formula in the case of polyhedral sets as it is recorded in \cite[Proposition 2.44]{primer}.

\begin{proposition}\label{pNormalsPolyhedral}{\rm (normal cone to polyhedral sets)}. Suppose that
  \[
  C = \{x\in \reals^n~|~Ax = b, ~Dx \leq d\},
  \]
  where $A$ and $D$ are $m\times n$ and $q\times n$ matrices, respectively. For any $\bar x\in C$,
\[
N_{C}(\bar x)=\big\{A^\top y + D^\top z~\big|~ y\in \reals^m; ~~z_i \geq 0 ~\mbox{ if } i\in \bbA(\bar x), ~~z_i = 0 \mbox{ otherwise}\big\},
\]
where $\bbA(\bar x) = \{i~|~\langle D_i, \bar x\rangle = d_i\}$, $D_i$ is the $i$th row of $D$, and $d_i$ is the $i$th component of $d$.
\end{proposition}

\begin{figure}
\centering
\includegraphics[width=0.61\textwidth]{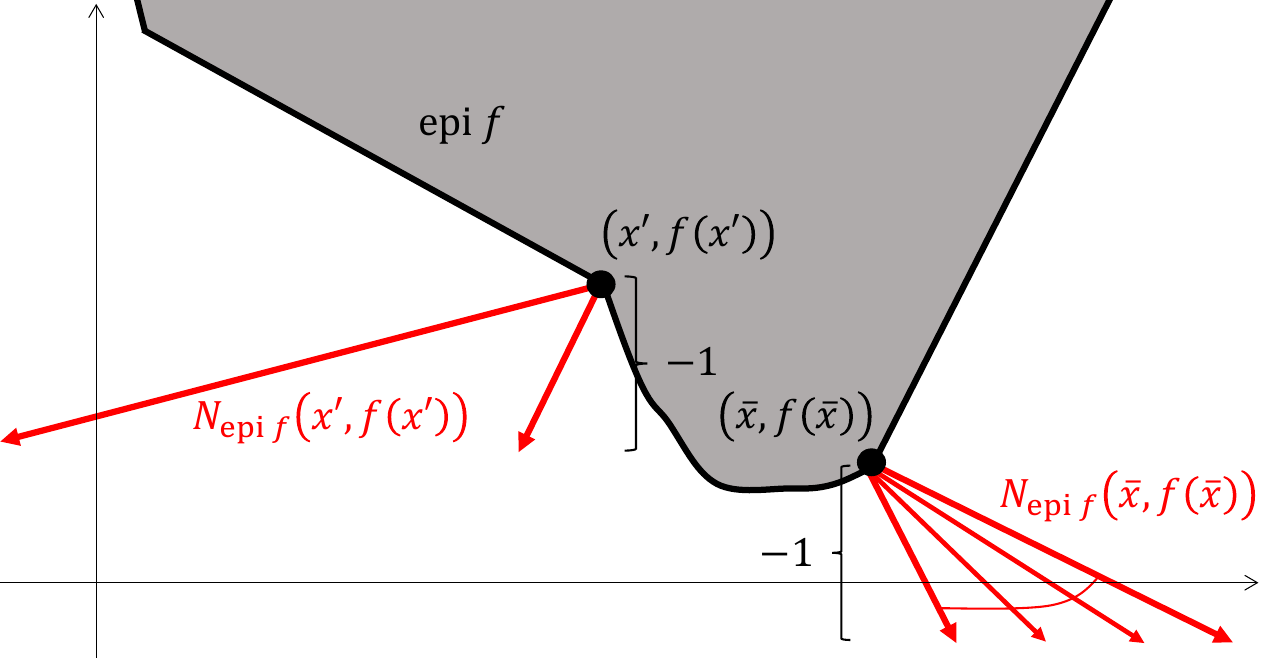}
\caption{Normal vectors of $\epi f$ defining the two subgradients $\partial f(x') = \{-4, -1/2\}$ and the interval of subgradients $\partial f(\bar x) = [1/2, 2]$.}\label{fig:subgrad}
\end{figure}

Subgradients now emerge from normal vectors of epigraphs. For a function $f:\reals^n\to \Reals$ and a point $\bar x$ where $f(\bar x)$ is finite, the {\em set of subgradients} of $f$ at $\bar x$ is 
\[
\partial f(\bar x) = \Big\{v\in \reals^n~\Big|~ (v,-1) \in  N_{\epi f}\big(\bar x, f(\bar x)\big)\Big\}.
\]
By convention, we set $\partial f(\bar x) = \emptyset$ when $f(\bar x)$ is not finite. A subgradient $v\in \partial f(\bar x)$ is also called a Mordukhovich, general, basic, or limiting subgradient. Figure \ref{fig:subgrad} shows the construction of subgradients of a nonsmooth nonconvex function. If $f$ is smooth in a neighborhood of $\bar x$, then $\partial f(\bar x)$ is a singleton with the gradient $\nabla f(\bar x)$ as its only element. If $f$ is convex and finite at $\bar x$, then $v\in \partial f(\bar x)$ if and only if $f(x) \geq f(\bar x) + \langle v, x-\bar x\rangle$ for all $x\in \reals^n$. Thus, the definition of subgradients extends classical notions from differential calculus and convex analysis. In the particular case of $h$, we find that (cf., for example, \cite[Section 5.I]{primer})
\begin{equation}\label{eqn:partialh}
\partial h(z) = \nargmin_{y\in Y} \tfrac{1}{2}\langle y, Qy\rangle - \langle y, z\rangle ~~\mbox{ for } z\in \dom h.
\end{equation}
A subgradient of $h$ is therefore computable using convex quadratic optimization, and it is unique if $Q$ is positive definite. An ability to compute subgradients easily is a motivation for the construction of $h$ in \eqref{eqn:h}.  We can also circle back to normal vectors and conclude that $\partial \iota_C(\bar x) = N_C(\bar x)$ for any $C\subset\reals^n$ and $\bar x\in \reals^n$.  

\begin{figure}
\centering
\includegraphics[width=0.51\textwidth]{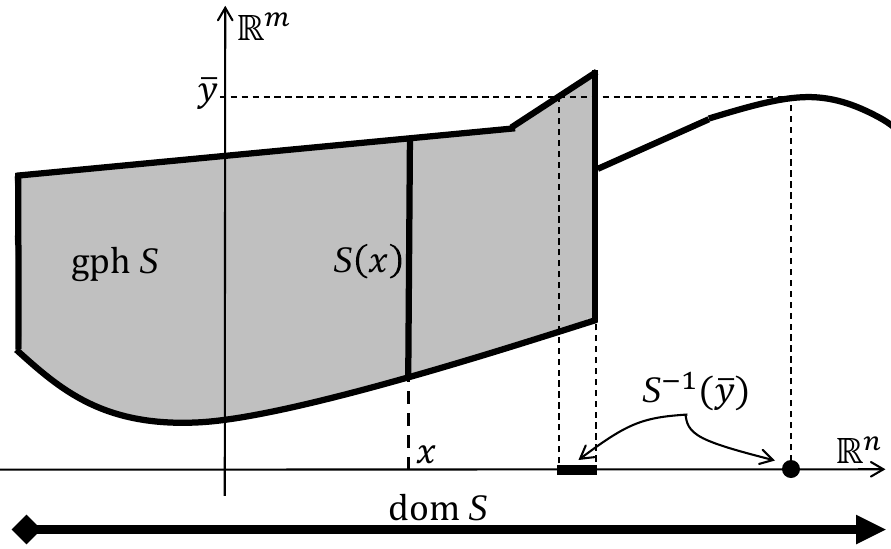}
\caption{The graph of a set-valued mapping $S:\reals^n\tto\reals^m$, its domain $\dom S$, and the solution set $S^{-1}(\bar y)$ to the generalized equation $\bar y \in S(x)$.}\label{fig:setvalued}
\end{figure}

We  introduce the concept of a {\em set-valued mapping}. While every argument of a function produces a scalar value (possibly $-\infty$ or $\infty$), a set-valued mapping $S:\reals^n\tto\reals^m$ yields at $x\in\reals^n$ the set $S(x)\subset\reals^m$. This multi-valuedness is indicated by the double arrows $\tto$. Figure \ref{fig:setvalued} illustrates the {\em graph} of $S$, which is defined by  $\gph S = \{(x,y) \in \reals^n \times \reals^m \, | \, y \in S(x) \}$, and also highlights the possibility that $S(x) = \emptyset$ or $S(x)$ is a singleton at some $x$. The {\em domain} of $S$ is $\dom S = \{x\in \reals^n \, | \, S(x) \neq \emptyset\}$. Examples of set-valued mappings include those defined by $S(x) = \partial f(x)$ and $S(x) = N_C(x)$. 

The set-valued mapping $S$ and a given point $\bar y\in \reals^m$ define a {\em generalized equation} $\bar y \in S(x)$. Its solution set
\[
S^{-1}(\bar y) = \big\{ x\in \reals^n ~\big|~ \bar y \in S(x) \big\}
\] 
is illustrated in Figure \ref{fig:setvalued}. The variational inequality in Example \ref{ex:VI} can be written as the generalized equation $0 \in S(x) = F(x) + N_C(x)$ because $C$ is convex and then $N_C(x)$ $=$ $\{v\in \reals^n ~|~\langle v, y - x\rangle \leq 0~\forall y\in C\}$; see, e.g., \cite[Proposition 4.42]{primer}.

We let $\nats = \{1, 2, \dots \}$ and typically index sequences by $\nu\in\nats$. We write $x^\nu \Nto x$ when a sequence $\{x^\nu, \nu\in\nats\}$ converges along a subsequence indexed by $\nu\in N\subset \nats$. 

The {\em point-to-set distance} from $\bar x\in\reals^n$ to $C\subset\reals^n$ is given by $\dist(\bar x, C) = \inf_{x \in C} \|x - \bar x\|_2$, with the convention that $\dist(\bar x, \emptyset) = \infty$.

\section{First-Order Optimality Conditions}\label{sec:firstorderoptimality}

A fundamental fact attributed to Pierre de Fermat \cite{Katz08} is that the derivative vanishes at every local minimizer of a smooth function. With the introduction of subgradients, an analogous property holds for arbitrary functions; see, e.g., \cite[Theorem 4.73]{primer} for an elementary proof.  

\begin{theorem}{\rm (Fermat rule).}\label{thm:Fermat}
For $f:\reals^n\to \Reals$ and $x^\star$ with $f(x^\star)\in \reals$, one has
\[
x^\star  \mbox{ is a local minimizer of } f ~\Longrightarrow~ 0\in \partial f(x^\star).
\]
\end{theorem}

While the Fermat rule only specifies a {\em necessary} condition for a local minimizer, it defines a generalized equation $0 \in \partial f(x)$ that can be solved as a substitute for the much more difficult task of finding a minimizer. It comes with the additional assurance that if $f$ is convex, then satisfying $0\in \partial f(x)$ is also sufficient for $x$ to be a minimizer.    

Practical use of Theorem \ref{thm:Fermat} requires calculus rules for computing subgradients. These are available in various forms; see, e.g., \cite[Chapter 10]{VaAn}, \cite[Chapter 4]{Mordukhovich.18}, and \cite[Section 4.I]{primer}. In the context of \eqref{eqn:actualproblem} and the stated assumptions, we obtain the following fact. 

\begin{theorem}{\rm (chain rule for composite function).}\label{thm:chainrule} For the function $f:\reals^n\to \Reals$ given by 
\[
f(x) = h\big(G(x)\big),
\]
with $h$ and $G$ as defined in Section \ref{sec:intro}, and a point $\bar x\in \dom f$, suppose that the following qualification holds:
\begin{equation}\label{eqn:compositeQual}
y\in N_{\dom h}\big(G(\bar x)\big) ~\mbox{ and }~ \nabla G(\bar x)^\top y = 0~~~\Longrightarrow~~~ y=0.
\end{equation}
Then, the set of subgradients of $f$ at $\bar x$ is given by 
\[
\partial f(\bar x) = \nabla G(\bar x)^\top \nargmin_{y\in Y} \tfrac{1}{2}\langle y, Qy\rangle - \big\langle y, G(\bar x)\big\rangle.
\]
\end{theorem}
\state Proof. Since $h$ is proper, lsc, and convex, one can invoke \cite[Theorem 4.64, Proposition 4.65]{primer} in conjunction with \eqref{eqn:partialh} to reach the conclusion.\eop

The qualification in \eqref{eqn:compositeQual} simplifies in many situations. If $Y$ is bounded as in Examples \ref{ex:goal},  \ref{ex:risk}, and \ref{ex:stat}, then $\dom h = \reals^m$ so that $N_{\dom h}(G(\bar x)) = \{0\}$. A positive definite $Q$ produces also a real-valued $h$. Thus, in either case, \eqref{eqn:compositeQual} holds automatically. If the $m\times n$ Jacobian matrix $\nabla G(\bar x)$ has rank $m$, then \eqref{eqn:compositeQual} also holds. Generally, 
\[
N_{\dom h}\big(G(\bar x)\big) = \big\{y \in Y^\infty ~\big|~Qy = 0, ~\big\langle y, G(\bar x)\big\rangle = 0\big\},
\]
where for an arbitrary $\bar y \in Y$, one has $Y^\infty = \{y\in \reals^m~|~\bar y + \lambda y \in Y ~\forall \lambda \in [0,\infty)\}$. Consequently, there are many avenues for verifying the qualification \eqref{eqn:compositeQual}. We refer to \cite{MohammadiMordukhovichSarabi.22} for a recent discussion of  \eqref{eqn:compositeQual}, its generalization to broader classes of compositions, and connections to {\em metric regularity} of epigraphs. 

The formula in Theorem \ref{thm:chainrule} shows that $\partial f(\bar x)$ is a singleton if and only if a convex quadratic optimization problem has a unique solution, which would be the case if $Q$ is positive definite. Regardless, all subgradients of $f$ at $\bar x$ are of the form $\nabla G(\bar x)^\top y$ for some $y\in Y$ determined by the convex quadratic optimization problem. 

It is possible to combine Theorem \ref{thm:chainrule} with a sum rule and thus incorporate $\iota_X$ for the purpose of addressing $\phi$ and developing a necessary optimality condition for \eqref{eqn:actualproblem2}; see, e.g., \cite[Theorem 4.75]{primer}. However, we proceed with a computationally motivated adjustment. When minimizing a smooth function $f$ using some algorithm, one might stop the algorithm if a current point $x$ satisfies $\|\nabla f(x)\|_2 \leq \epsilon$, where $\epsilon$ is a positive tolerance. The same strategy is problematic when $f$ is nonsmooth. For example, if $f(x) = |x|$, then $\partial f(x) = \{-1\}$ when $x<0$, $\partial f(x) = [-1,1]$ when $x=0$, and $\partial f(x) = \{1\}$ otherwise. In this case, $\partial f(x)$ is bounded away from the origin unless $x=0$. Thus, the condition $\dist(0, \partial f(x))\leq \epsilon$ will not be suitable as a stopping criterion. The following optimality condition addresses such concerns.

\begin{theorem}{\rm (optimality condition for actual problem).}\label{thm:OptimComposite}
Suppose that $\bar x$ is a local minimizer of \eqref{eqn:actualproblem}, with $X$, $Y$, $Q$, $h$, and $G$ as defined in Section \ref{sec:intro}, and that the following qualification holds:
\begin{equation}\label{eqn:compositeQual2}
y\in N_{\dom h}\big(G(\bar x)\big) ~\mbox{ and }~ -\nabla G(\bar x)^\top y \in N_X(\bar x)~~~\Longrightarrow~~~ y=0.
\end{equation}
Then, there are $\bar y\in \reals^m$ and $\bar z\in \reals^m$ such that 
\[
0 \in \Phi(\bar x,\bar y,\bar z),
\]
where the set-valued mapping $\Phi:\reals^{n+2m}\tto\reals^{2m +n}$ is given by 
\begin{equation}\label{eqn:Phi}
\Phi(x,y,z) = \big\{G(x) - z\big\} \times \big(Qy - z + N_Y(y) \big) \times \big(\nabla G(x)^\top y + N_X(x)\big).
\end{equation}
\end{theorem}
\state Proof. By \cite[Theorem 4.75]{primer}, there exists $\bar y\in \partial h(G(\bar x))$ such that $-\nabla G(\bar x)^\top \bar y \in N_X(\bar x)$. Using \eqref{eqn:partialh}, we conclude that $\bar y\in \partial h(G(\bar x))$ if and only if $G(\bar x) - Q\bar y \in N_Y(\bar y)$. With $\bar z = G(\bar x)$, the conclusion follows.\eop

The auxiliary vectors $y$ and $z$ in the optimality condition introduce additional flexibility that is algorithmically important. The problematic case of $|\cdot|$ discussed prior to the theorem is now handled naturally. Let $X = \reals$, $h(z) = |z| = \sup_{y\in [-1,1]} yz$, and $G(x) = x$. Then, $\phi(x) = \iota_X(x) + h(G(x)) = |x|$ and the optimality condition of Theorem \ref{thm:OptimComposite} becomes
\[
0 \in \Phi(x,y,z) = \{x - z\} \times \big(\hspace{-0.06cm}-z + N_{[-1,1]}(y)\big) \times \{y\} ~\mbox{ for some $y$ and $z$.  }
\] 
Now, the stopping criterion $\dist(0, \Phi(x,y,z)) \leq \epsilon$ is meaningful. It would be satisfied for $x$ near $0$ and not only at $x = 0$. Additional motivation for the optimality condition of Theorem \ref{thm:OptimComposite} emerges in the discussion of data inaccuracies in Section \ref{sec:approx}.    

There are numerous algorithms that can solve the generalized equation $0 \in \Phi(x,y,z)$; see, e.g., \cite[Chapter 7]{primer}. If not only $Y$ but also $X$ is a polyhedral set, then Proposition \ref{pNormalsPolyhedral} furnishes convenient expressions for both $N_X(x)$ and $N_Y(y)$ and solvers such as PATH \cite{DirkseFerris.95} become available when $G$ is twice smooth (i.e., of class $C^2$).

The discussion after Theorem \ref{thm:chainrule} provides guidance about when the qualification \eqref{eqn:compositeQual2} holds. In the setting of Example \ref{ex:nonlinear}, the qualification \eqref{eqn:compositeQual2} is equivalent to the Mangasarian--Fromovitz constraint qualification. We refer to \cite{MohammadiMordukhovichSarabi.22} for an up-to-date treatment of qualifications in the setting of composite optimization and to \cite{DoMordukhovichSarabi.21} for a discussion about criticality of multiplier vectors.

\section{Approximations and Consistency}\label{sec:approx}

Most real-world instances of the actual problem \eqref{eqn:actualproblem2} are subject to data ambiguity. The target values $\tau_1, \dots, \tau_m$ in Example \ref{ex:goal} or the risk-averseness parameter $\alpha$ in Example \ref{ex:risk} might not be fully known. In other situations, we introduce approximations to facilitate computations. Examples \ref{ex:goal}-\ref{ex:VI} have $Q$ in \eqref{eqn:h} as a zero matrix, but it could be approximated by positive definite matrices $\{Q^\nu = I/\nu, \nu\in\nats\}$, where $I$ is the $m\times m$ identity matrix. This change has the benefit that the resulting approximations $h^\nu(z) = \sup_{y\in Y} ~\langle y, z\rangle - 1/2\langle y, Q^\nu y\rangle$ are smooth as can be seen from \eqref{eqn:partialh}. 

Regardless of the cause, what is the effect of approximating data on problem solutions? In the setting of Example \ref{ex:nonlinear} with $g_0(x) = -x$, a single inequality constraint $g_1(x) = x^3 - x^2 - x + 1\leq 0$, and no equality constraints, the unique minimizer becomes $x^\star = 1$. This shifts to near $-1$ after the change from $g_1$ to $g_1^\nu(x) =  x^3 - x^2 - x + 1 + 1/\nu$. While $g_1^\nu$ converges uniformly to $g_1$ as $\nu\to \infty$, the corresponding minimizers do not converge. One might turn to a representation of constraints using an indicator function as in \eqref{eqn:actualproblem2} and consider uniform convergence or pointwise convergence but this is also problematic; see \cite[Section 4.C]{primer} for details. In summary, these classical notions of convergence are not suitable for analysis of optimization problems and generalized equations. Variational analysis brings in another concept.   

A sequence of sets $\{C^\nu\subset\reals^n, \nu\in\nats\}$ converges to $C\subset\reals^n$ in the sense of {\em Painlev\'{e}--Kuratowski} when $C$ is closed and $\dist(x,C^\nu)\to \dist(x,C)$ for all $x\in\reals^n$. It is denoted by $C^\nu\sto C$ and referred to as {\em set-convergence}. 

More specifically, we say that a sequence of functions $\{f^\nu:\reals^n\to \Reals, \nu\in\nats\}$ {\em epi-converges} to $f:\reals^n\to \Reals$ if $\epi f^\nu \sto \epi f$. It turns out that epi-convergence is the more natural concept to consider when analyzing approximations of minimization problems. The following facts are compiled from \cite[Theorem 5.5]{primer} and \cite[Proposition 7.7]{VaAn}; see the commentary in \cite[Chapter 7]{VaAn} for their historical development.

\begin{theorem}{\rm (consequences of epi-convergence).}\label{thm:epiconv}
For $f,f^\nu:\reals^n\to \Reals$, $\nu\in\nats$, if $\epi f^\nu \sto \epi f$, then the following hold. 
\begin{enumerate}[(a)]

\item Suppose that $\inf f < \infty$,  $\epsilon^\nu\to 0$, and $x^\nu \in \epsilon^\nu\mbox{-}\nargmin f^\nu$ for $\nu\in\nats$. If there is a subsequence $N\subset \nats$ such that $x^\nu\Nto x$, then $x\in\nargmin f$ and $f^\nu(x^\nu) \Nto \inf f$. 
    
\item Suppose that $\alpha^\nu\to \alpha\in\reals$ and $x^\nu$ satisfies $f^\nu(x^\nu) \leq \alpha^\nu$ for $\nu\in\nats$. If there is a subsequence $N\subset \nats$ such that $x^\nu\Nto x$, then $f(x) \leq \alpha$. 

\item For any $x\in \reals^n$ and $\alpha\in\reals$ satisfying $f(x) \leq \alpha$, there exist $x^\nu\to x$ and $\alpha^\nu\to \alpha$ such that $f^\nu(x^\nu) \leq \alpha^\nu$ for all $\nu$. 

\end{enumerate}    
\end{theorem}    
    
A minimization problem represented by a function $f$ is suitably approximated by other problems expressed in terms of $f^\nu$ when $\epi f^\nu \sto \epi f$. Not only are minimizers well-behaved as in Theorem \ref{thm:epiconv}(a), but level-sets of $f^\nu$ are close to those of $f$ as seen from parts (b), (c). Thus, epi-convergence entails a notion of ``global'' approximation that goes much beyond the often cited fact in part (a). It is unencumbered by extended real-valuedness, nonconvexity, or nonsmoothness of the functions involved. 

A set-valued mapping $S:\reals^n\tto\reals^m$, for example given by subgradient and/or normal cone mappings, may also be approximated due to data ambiguity and computational considerations. What is the effect on the solutions of the resulting generalized equations? Again, set-convergence enters as a key concept but now applied to the graphs of set-valued mappings, which we referred to as {\em graphical convergence}. The following consequences of graphical convergence can be found, e.g, in \cite[Theorem 7.34]{primer}.   
    
\begin{theorem}{\rm (consequences of graphical convergence).}\label{thm:graphconv}
For $S,S^\nu:\reals^n\tto \reals^m$, $\nu\in\nats$, suppose that $\gph S^\nu \sto \gph S$. Then, the following hold. 
\begin{enumerate}[(a)]

\item If $x^\nu \Nto x$ for some subsequence $N\subset \nats$ and $\dist(0,S^\nu(x^\nu)) \leq \epsilon^\nu$, then $0 \in S(x)$ whenever $\epsilon^\nu\to 0$. 
    
\item If $x$ satisfies $0 \in S(x)$, then there exist $\epsilon^\nu\to 0$ and $x^\nu\to x$ with $\dist(0, S^\nu(x^\nu)) \leq \epsilon^\nu$ for all $\nu$. 

\end{enumerate}    
\end{theorem}    

The theorem establishes that any approximations that might be introduced to a generalized equation $0\in S(x)$ have vanishing effect on the solutions when they can be represented by set-valued mappings $S^\nu$ graphically converging to $S$. 

In the context of the actual problem \eqref{eqn:actualproblem2}, epi-convergence and graphical convergence provide a comprehensive framework for studying approximations regardless of their origin. For each $\nu\in\nats$, suppose that $X^\nu\subset\reals^n$ is a nonempty, closed, and convex set approximating $X$, $G^\nu:\reals^n\tto\reals^m$ is a smooth mapping approximating $G$, $Q^\nu$ is an $m\times m$ symmetric positive semidefinite matrix replacing $Q$, and $Y^\nu\subset\reals^m$ is a nonempty polyhedral set replacing $Y$. Then, the {\em approximating problem}  
\begin{equation}\label{eqn:approxproblem}
\nnmin_{x\in \reals^n} \,\phi^\nu(x)= \iota_{X^\nu}(x) + h^\nu\big(G^\nu(x)\big), ~~\mbox{ where } h^\nu(z) = \sup_{y\in Y^\nu} \langle y,z\rangle - \tfrac{1}{2} \langle y, Q^\nu y\rangle,
\end{equation}
is of the same form as \eqref{eqn:actualproblem2} but with approximating data. Just as $0\in \Phi(x,y,z)$ furnishes a necessary optimality condition for \eqref{eqn:actualproblem2} under a qualification (cf. Theorem \ref{thm:OptimComposite}), we define 
\begin{equation}\label{eqn:Phinu}
\Phi^\nu(x,y,z) = \big\{G^\nu(x) - z\big\} \times \big(Q^\nu y - z + N_{Y^\nu}(y) \big) \times \big(\nabla G^\nu(x)^\top y + N_{X^\nu}(x)\big)
\end{equation}
so that $0\in \Phi^\nu(x,y,z)$ is a parallel necessary condition for the approximating problem \eqref{eqn:approxproblem}.  

Following \cite{Royset.22a} and \cite[Section 7.I]{primer}, we say that the pairs $\{(\phi^\nu, \Phi^\nu), \nu\in \nats\}$ from \eqref{eqn:approxproblem} and \eqref{eqn:Phinu} are {\em consistent approximations} of $(\phi,\Phi)$ from \eqref{eqn:actualproblem2} and \eqref{eqn:Phi} when
\[
\epi \phi^\nu \sto \epi \phi~~~~\mbox{ and }~~~~\gph \Phi^\nu \sto \gph \Phi.
\]

It is apparent that under consistency the approximating problems \eqref{eqn:approxproblem} are ``good'' approximations of the actual problem \eqref{eqn:actualproblem2}. Minimizers, minimum values, level-sets, and points satisfying optimality conditions are all well-behaved under approximations in the sense of Theorems \ref{thm:epiconv} and \ref{thm:graphconv}. Algorithmically, one might leverage consistency by solving a sequence of approximating problems instead of tackling the actual problem directly. 

\medskip

\state Consistent Approximation Algorithm.

\begin{description}

  \item[Data.] ~\,~$\epsilon^\nu \geq 0$, with $\epsilon^\nu\to 0$.

  \item[Step 0.]  Set $\nu = 1$.

  \item[Step 1.]  Apply an algorithm to \eqref{eqn:approxproblem} until it obtains $x^\nu$, with corresponding $(y^\nu,z^\nu)$, satisfying
  \[
  \dist\big(0, \Phi^\nu(x^\nu,y^\nu,z^\nu)\big) \leq \epsilon^\nu.
\]

\item[Step 2.] Replace $\nu$ by $\nu +1$ and go to Step 1.
\end{description}

\medskip

With a judicious choice of approximating problems, Step 1 might be accomplished efficiently using existing algorithms. Any cluster point $(x,y,z)$ of the resulting sequence $\{(x^\nu,y^\nu,z^\nu), \nu\in\nats\}$ would then satisfy $0\in \Phi(x,y,z)$ provided that $\gph \Phi^\nu \sto \gph \Phi$; see Theorem \ref{thm:graphconv}. Consistency strengthens this conclusion further by guaranteeing properties for minimizers, minimum values, and level-sets via Theorem \ref{thm:epiconv}. 

The objective functions $\phi^\nu$ of the approximating problems \eqref{eqn:approxproblem} and their optimality conditions expressed by $\Phi^\nu$ in \eqref{eqn:Phinu} are consistent under natural assumptions.

\begin{theorem}{\rm (consistency).}\label{thm:consistency} In the setting of \eqref{eqn:actualproblem2} and \eqref{eqn:approxproblem}, suppose that $X^\nu\sto X$, $Y^\nu\sto Y$, $Q^\nu\to Q$, and $G^\nu(x^\nu)\to G(x)$ as well as $\nabla G^\nu(x^\nu)\to \nabla G(x)$ whenever $x^\nu\to x$. Then, $\gph \Phi^\nu \sto \gph \Phi$. 

Additionally, $\epi \phi^\nu\sto \epi\phi$ provided that one of the following conditions also holds:  
\begin{enumerate}[(a)]

\item $Y$ is bounded, 
    
\item $Q$ is positive definite, 
    
\item $Y^\nu \subset Y$, $X^\nu = X$, and $G^\nu = G$ for all $\nu\in\nats$.
    
\end{enumerate}    

Consequently, $\{(\phi^\nu, \Phi^\nu), \nu\in \nats\}$ are {\em consistent approximations} of $(\phi,\Phi)$ if (a), (b), or (c) holds. 

\end{theorem}    
\state Proof. The conclusion about $\gph \Phi^\nu \sto \gph \Phi$ follows from \cite[Theorem 7.45(c)]{primer} after realizing that $\epi h^\nu\sto \epi h$ under the stated assumptions; see \cite[Exercise 7.47]{primer}. Under (a) and (b), $\epi \phi^\nu\sto \epi\phi$ follows by \cite[Theorem 7.45(a)]{primer}. Under (c), $Y^\nu \subset Y$ implies that $h^\nu(z) \leq h(z)$ for all $z$. This fact together with the general consequence $\nliminf h^\nu(z) \geq h(z)$ from  epi-convergence (see, e.g., \cite[Theorem 4.15(a)]{primer}) imply that $h^\nu$ converges pointwise to $h$. We can then bring in \cite[Theorem 7.45(b)]{primer} to reach the conclusion.\eop  

Consistency for various approximations follows immediately in Examples \ref{ex:goal} and \ref{ex:risk} using Theorem \ref{thm:consistency}(a) because $Y$ is bounded. In the context of Example \ref{ex:nonlinear}, an approximation that falls under Theorem \ref{thm:consistency}(c) is to replace 
\[
Y = \{1\} \times \reals^m \times [0, \infty)^q ~~\mbox{ by }  ~~Y^\nu = \{1\} \times [-\theta^\nu, \theta^\nu]^m \times [0, \theta^\nu]^q
\]
for some $\theta^\nu \in [0,\infty)$, Then, with $y = (y_0, y_1, \dots, y_{m+q})$ and $z = (z_0, z_1, \dots, z_{m+q})$,
\[
h^\nu(z) = \sup_{y\in Y^\nu} \langle y, z\rangle = z_0 + \theta^\nu \sum_{i=1}^m |z_i| + \theta^\nu \sum_{i=m+1}^{m+q} \max\{0, z_{i}\},  
\]
which results in the classical exact penalization of equality and inequality constraints. If $\theta^\nu\to \infty$, then $Y^\nu\sto Y$ and consistency follows by Theorem \ref{thm:consistency}(c). We refer to \cite{Royset.22a} for additional examples and refinements. 

Approximations of the matrix $Q$ by the positive definite matrices $\{Q^\nu = Q + I/\nu, \nu\in\nats\}$ hinted to in the introductory paragraph of this section satisfy $Q^\nu\to Q$ and produce $h^\nu(z) = \sup_{y\in Y} ~\langle y, z\rangle - 1/2\langle y, Q^\nu y\rangle$. At least when $Y$ is bounded or $Q$ is positive definite, the function $h^\nu$ turns out to be a {\em Moreau envelope} of $h$, a general functional approximation tool frequently appearing in variational analysis; see, e.g, \cite[Section 1.G]{VaAn} and \cite{AttouchWets.91,BurkeHoheisel.13,LiCui.24}. For details about the specific relationship between $h$ and $h^\nu$ viewed through the lens of Moreau envelopes, we refer to \cite[Proposition 4.11]{BurkeHoheisel.13}.

\section{Algorithms}\label{sec:algo}

With the vast number of applications covered by the actual problem \eqref{eqn:actualproblem2}, one cannot expect a single algorithm to address all its instances effectively. We convey three algorithmic ideas, with additional possibilities emerging in the following sections.\\  

\noindent {\bf Epigraphical reformulation.} For specific $Y$ and $Q$, with $G$ being an affine mapping, there is a reformulation that take us back to linear, quadratic, or convex optimization via the epigraphical reformulation: 
\begin{equation}\label{eqn:epigraphreform}
h\big(G(x)\big) \leq \alpha ~~~\Longleftrightarrow~~~\big\langle y,G(x)\big\rangle - \tfrac{1}{2}\langle y, Qy\rangle \leq \alpha ~~\forall y\in Y.
\end{equation}
An approach for solving \eqref{eqn:actualproblem2} would be to minimize $\alpha$ subject to these constraints as well as $x\in X$ using an existing solver. This approach is viable if the infinite collection of constraints indexed by $y\in Y$ can be expressed using a moderate number of constraints, at least in an approximate sense. The reformulation \eqref{eqn:epigraphreform} holds also when $G$ is not affine, but then the approach needs to contend with a potentially large number of nonlinear inequality constraints and this could cause difficulty for solvers.\\

\noindent {\bf Dual reformulation.} Another approach to solving \eqref{eqn:actualproblem2} leverages an alternative expression for $h$ developed in \cite{Rockafellar.99} (see also \cite[Section 5.I]{primer}) by using quadratic programming duality. Suppose that 
\begin{equation}\label{eqn:Yform}
Y = \{y \in \reals^m~|~A^\top y \leq b\} ~\mbox{ and } ~Q = D J^{-1} D^\top
\end{equation}
for some $m\times q$-matrix $A$, vector $b\in \reals^q$, $m\times m$-matrix $D$, and symmetric positive definite $m\times m$-matrix $J$. Since every polyhedral set can be written in terms of a finite collection of linear inequality constraints, the requirement on $Y$ is not restrictive. The requirement on $Q$ is also reasonable. If $Q=0$, we can select $D=0$ and $J=I$, the $m\times m$ identity matrix. If $Q$ is positive definite, then spectral decomposition gives $D$ and $J$. Specifically, there are an $m\times m$-matrix $D$ (consisting of orthonormal eigenvectors) and an $m\times m$-matrix $\Lambda$, with the eigenvalues of $Q$ along its diagonal and zero elsewhere, such that $Q = D \Lambda D^\top$. Since $Q$ is symmetric and positive definite, its eigenvalues are positive and we can set $J=\Lambda^{-1}$.

With these formulas for $Y$ and $Q$, quadratic programming duality produces the formula
\begin{equation}\label{eqn:hreform}
h(z)= \inf_{v\in \reals^q, w\in \reals^m} \Big\{\langle b,v\rangle + \tfrac{1}{2} \langle w, Jw\rangle ~\Big|~ A v + Dw = z, ~ v \geq 0\Big\}.
\end{equation}
Thus, \eqref{eqn:actualproblem2} is equivalent to solving 
\[
\nnmin_{x\in X, v\in \reals^q, w\in \reals^m}  \langle b,v\rangle + \tfrac{1}{2} \langle w, Jw\rangle + \tfrac{1}{2\lambda^\nu} \|x-x^\nu\|_2^2
~~\mbox{ subject to }~ ~  A v + Dw = G(x), ~~ v \geq 0,
\]
which then involves auxiliary vectors $v$ and $w$ to be discarded after optimization. When $G$ is highly nonlinear, however, the resulting $m$ nonlinear equality constraints $A v + Dw = G(x)$ may be numerically challenging to handle by standard nonlinear optimization solvers.

\begin{example}{\rm (risk minimization (cont.)).}\label{ex:risk2} Returning to Example \ref{ex:risk}, we find that $Y = \{y \in \reals^m~|~A^\top y \leq b\}$ with $b = (1, -1, 0, \dots, 0, p_1/(1-\alpha), \dots, p_m/(1-\alpha))\in \reals^{2+2m}$ and 
\[
A = \begin{bmatrix}
  1      & -1     & -1     & 0      & \cdots & 0      & 1      & 0      & \cdots & 0\\ 
  1      & -1     & 0      & -1     & \cdots & 0      & 0      & 1      & \cdots & 0\\ 
  \vdots & \vdots & \vdots &        & \ddots & \vdots & \vdots &        & \ddots & \vdots\\ 
  1      & -1     & 0      & \cdots & 0      &-1      & 0      & \cdots & 0      & 1
\end{bmatrix}.
\]
Since $Q=0$ in this example, one can set $D = 0$ and $J$ equal to the identity matrix in the reformulation \eqref{eqn:hreform}. Thus, the optimal $w$ in \eqref{eqn:hreform} is necessarily the zero vector. This simplification and some other ones afforded by the particular $A$-matrix and $b$-vector yield that   
\[
h(z)= \inf_{\gamma\in\reals, u\in \reals^m} \Big\{ \gamma + \frac{1}{1-\alpha}\sum_{i=1}^m p_i u_i ~\Big|~ z_i - \gamma \leq u_i, ~u_i\geq 0, ~~i=1, \dots, m\Big\}.
\]
Thus, the minimization of $h(G(x))$ over $x\in X$ in this example can be accomplished by solving 
\[
\nnmin_{x\in X, \gamma\in\reals, u\in \reals^m} \gamma + \frac{1}{1-\alpha}\sum_{i=1}^m p_i u_i ~\mbox{ subject to } ~ g_i(x) - \gamma \leq u_i, ~u_i\geq 0, ~~i=1, \dots, m.
\]
This formulation is identical to a common one emerging from the Rockafellar-Uryasev formula for $\alpha$-superquantiles in \cite{RockafellarUryasev.00}.
\end{example} 

\medskip

\noindent {\bf Proximal composite method.} The classical gradient descent method of Cauchy motivates another approach. To minimize a smooth function $f:\reals^n\to \reals$, one might start with $x^0\in\reals^n$ and then iterate using 
\[
x^{\nu+1} = x^\nu - \lambda^\nu \nabla f^\nu(x^\nu) ~~\Longleftrightarrow~~ x^{\nu+1} \in \nargmin_x f(x^\nu) + \big\langle \nabla f(x^\nu), x - x^\nu\big\rangle + \tfrac{1}{2\lambda^\nu} \|x - x^\nu\|_2^2, 
\]
where $\lambda^\nu$ is a positive step size. The right-most interpretation of the gradient descent step highlights two characteristics: the function $f$ is linearized at the current point and a quadratic penalty prevents any step from being too long, which might jeopardize the accuracy of the linear approximation of the function. We can leverage these ideas when addressing \eqref{eqn:actualproblem2}: an algorithm emerges from linearizing $G$ and adding a quadratic penalty term.   

\bigskip

\state Proximal Composite Method.

\begin{description}

  \item[Data.] ~\,~$x^0\in X$, $\tau\in (1,\infty)$, $\sigma\in (0,1)$, $\bar \lambda \in (0,\infty)$, $\lambda^0 \in (0, \bar \lambda]$.

  \item[Step 0.]  Set $\nu = 0$.

  \item[Step 1.]  Compute
\[
  \bar x^\nu \in \nargmin_{x\in X} h\big(G(x^\nu) + \nabla G(x^\nu)(x-x^\nu)\big) + \tfrac{1}{2\lambda^\nu}\|x-x^\nu\|_2^2.
\]
~~~~~If $\bar x^{\nu}=x^\nu$, then Stop.

  \item[Step 2.]  If
  \[
  h\big(G(x^\nu)  \big) - h\big(G(\bar x^\nu)  \big) \geq \sigma \Big( h\big(G(x^\nu)\big) - h\big(G(x^\nu) + \nabla G(x^\nu)(\bar x^\nu - x^\nu) \big) \Big),
  \]
~~~~~then set $\lambda^{\nu+1} = \min\{\tau \lambda^\nu, \bar\lambda\}$ and go to Step 3.

~~~~~Else, replace $\lambda^\nu$ by $\lambda^\nu/\tau$ and go to Step 1.

  \item[Step 3.] Set $x^{\nu+1} = \bar x^\nu$, replace $\nu$ by $\nu +1$ and go to  Step 1.
\end{description}

\smallskip

Every cluster point produced by the proximal composite method satisfies the optimality condition of Theorem \ref{thm:OptimComposite} under two additional assumptions: The mapping $G$ is actually twice smooth and $h$ is real-valued. This fact is summarized next.

\begin{theorem}{\rm (proximal composite method).}\label{cProxComposite} In the setting of the actual problem \eqref{eqn:actualproblem2} and its optimality condition given by $\Phi$ from \eqref{eqn:Phi}, suppose that $G$ is twice smooth, $h$ is real-valued, and the proximal composite method has generated a sequence $\{x^\nu, \nu\in \nats\}$ with a cluster point $x^\star$. Then, there are $y^\star \in \reals^m$ and $z^\star\in \reals^m$ such that $0\in \Phi(x^\star,y^\star,z^\star)$.
\end{theorem}
\state Proof. By convergence statement 6.35 in \cite{primer}, which in turn specializes \cite{LewisWright.16}, there is $y^\star \in \partial h(G(x^\star))$ such that $-\nabla G(x^\star)^\top y^\star \in N_X(x^\star)$. The result follows by arguing as in the proof of Theorem \ref{thm:OptimComposite}.\eop

Since the approximating problems \eqref{eqn:approxproblem} are of the same form as the actual problem \eqref{eqn:actualproblem2}, the proximal composite method could also be applied to \eqref{eqn:approxproblem}. This makes the additional requirements of twice smoothness and real-valuedness less restrictive as such properties might be {\em constructed} in the approximating problems by smoothing $\nabla G$ and bounding $h$. One can then combine the consistent approximation algorithm of Section \ref{sec:approx} with the proximal composite method. Step 1 of the consistent approximation algorithm would amount to several iterations of the proximal composite method applied to a specific approximating problem; see \cite{Royset.22a} for further details. 

Regardless of the context, the subproblem in Step 1 of the proximal composite method is convex. The best way of solving it depends on further details about $h$. If the reformulation \eqref{eqn:epigraphreform} is viable, then it furnishes an approach to accomplish Step 1 using standard convex optimization algorithms or QP-solvers if $X$ is polyhedral. If $Y$ is in the form \eqref{eqn:Yform}, then the reformulation \eqref{eqn:hreform} becomes available and Step 1 amounts to solving 
\begin{align*}
\nnmin_{x\in X, v\in \reals^q, w\in \reals^m} & \langle b,v\rangle + \tfrac{1}{2} \langle w, Jw\rangle + \tfrac{1}{2\lambda^\nu} \|x-x^\nu\|_2^2\\
\mbox{ subject to } ~ & A v + Dw = G(x^\nu) + \nabla G(x^\nu)(x-x^\nu), ~~ v \geq 0.
\end{align*}
Again, one can leverage convex optimization algorithms or QP-solvers if $X$ is polyhedral. 

Proximal composite methods can be traced back to \cite{Fletcher.82,Powell.83,Burke.85}, with trust-region versions appearing in \cite{Powell.84,Yuan.85,BurkeFerris.95}. More recent refinements include \cite{Sagastizabal.13} for positively homogeneous real-valued $h$ (as in Examples \ref{ex:goal}, \ref{ex:risk}, and \ref{ex:stat}), \cite{LewisWright.16} to handle extended real-valued $h$, and \cite{DrusvyatskiyPaquette.19} to address rate of convergence and inexact solution of subproblems.

\section{Problem Relaxation and Duality}\label{sec:duality}

An optimization problem can be associated with many substitute problems that offer new insights and computational opportunities. Among these the dual problems of linear and conic optimization are utilized in analysis and algorithms. We adopt a much broader point of view that includes these possibilities as special cases but extends further to nonconvex problems. 

The key concept (pioneered in \cite[Chapter 29]{Rockafellar.70} and \cite{Rockafellar.74,Rockafellar.85}) is that of problem parametrization as encoded by a Rockafellian function. Following the terminology in \cite[Section 5.A]{primer}, we say that $f:\reals^q\times\reals^n\to \Reals$ is a {\em Rockafellian} for a given function $g:\reals^n\to\Reals$ when $f(0,x) = g(x)$ for all $x\in\reals^n$. 

A Rockafellian $f$ allows us to supplement the minimization of $g$ by alternative problems constructed from $f$, with minimizing $f(0, \cdot)$ bringing us back to the original problem. We summarize the main possibilities: 
\begin{align*}
\mbox{(original problem)}~~ & ~~\nnmin_{x\in \reals^n} \, g(x)~~\Longleftrightarrow~~\nnmin_{x\in\reals^n} \, f(0,x)\\
\mbox{(Rockafellian relaxation)}~~ & ~~\nnmin_{u\in\reals^q, x\in \reals^n} \, f(u,x) - \langle y, u\rangle\\
\mbox{(Lagrangian relaxation)}~~ & ~~\nnmin_{x\in \reals^n} \, l(x,y) = \inf_{u\in\reals^q} f(u,x) - \langle y, u\rangle\\
\mbox{(dual problem)}~~ &~~ \nnmax_{y\in \reals^q} \, \psi(y) = \inf_{x\in\reals^n} l(x,y).
\end{align*}
Rockafellian relaxation and Lagrangian relaxation are indeed relaxations of the original problem in the sense that their minimum values are no greater than that of the original problem. Moreover, they rely on a {\em multiplier vector} $y \in \reals^q$, which in turn is optimized in the dual problem. A dual problem therefore defines the best lower bound of $\inf g$ under the adopted Rockafellian. In this section, we apply these general principles in the context of the actual problem from Section \ref{sec:intro}.

As for any optimization problem, the actual problem \eqref{eqn:actualproblem2} can be associated with many Rockafellians. We consider two possibilities.  

The first possibility is the function $f:\reals^m\times \reals^n\to \Reals$ given by 
\begin{equation}\label{eqn:Rock}
f(u,x) = \iota_X(x) + h\big(G(x) + u\big),  
\end{equation}
which indeed is a Rockafellian for $\phi$ in \eqref{eqn:actualproblem2} because $f(0,x) = \phi(x)$ for all $x\in \reals^n$. An advantage of this Rockafellian is that the corresponding Lagrangian $l:\reals^n\times\reals^m\to \Reals$ has the specific form
\begin{equation}\label{eqn:lagr}
l(x,y) = \iota_X(x) + \big\langle G(x),y\big\rangle - \tfrac{1}{2} \langle y, Qy\rangle - \iota_Y(y)
\end{equation}
because the minimization of $f(u,x) - \langle y, u\rangle$ over $u$ can be carried out explicitly; see, e.g., \cite[Proposition 5.28, Section 5.I]{primer}. 
The formula for $l(x,y)$ may involve the unorthodox arithmetic operation ``$\infty-\infty$'' and this requires clarification. In the context of minimization (as in minimizing a Lagrangian with respect to $x$), variational analysis adopts the convention $\infty - \infty = \infty$. Thus, $l(x,y) = \infty$ for every $x\not\in X$ regardless of $y$. We refer to \cite[p.13]{primer} for further details.  
  
By definition, a simplified formula for the dual objective function $\psi$ stems from an explicit calculation of the minimum value of the Lagrangian with respect to $x$. Naturally, this is rarely available. If $G(x) = b - Ax$, however, then we obtain 
\[
\psi(y) = -\iota_Y(y) + \langle b, y\rangle - \tfrac{1}{2} \langle y, Qy\rangle - \tilde h(A^\top y), ~\mbox{ where } \tilde h(v) = \sup_{x\in X} \langle v,x\rangle. 
\]   
Interestingly, $\tilde h$ is nearly of the same form as $h$ in \eqref{eqn:h}. A more complete symmetrical relationship between minimizing $\phi$ and maximizing $\psi$ is in fact available; see, e.g., \cite[Example 5.57]{primer}.

While Rockafellian relaxations, Lagrangian relaxations, and dual problems always furnish lower bounds on the minimum value of the actual problem, these bounds are tight under additional assumptions. 

\begin{theorem}{\rm (strong duality).}\label{thm:duality}  
Under the Rockafellian in \eqref{eqn:Rock} for $\phi$ in \eqref{eqn:actualproblem2}, the corresponding dual function $\psi:\reals^m\to \Reals$ satisfies  
\[
\sup_{y\in \reals^m} \psi(y) = \inf_{x\in \reals^n} \phi(x)
\] 
provided that any one of the following conditions holds: 
\begin{enumerate}[(a)]
  \item $\langle G(\cdot), y\rangle$ is convex for all $y \in Y$ and, in addition, either $X$ or $Y$ is compact;

  \item $X$ is polyhedral, $G(x) = b - Ax$, and, in addition, either $\inf_{x\in \reals^n} \phi(x)$ or $\sup_{y\in \reals^m} \psi(y)$ is finite.
\end{enumerate}
\end{theorem}
\state Proof. Part (a) follows by Sion’s theorem \cite[Corollary 3.3]{Sion.58}; (b) holds by \cite[Theorem 11.42]{VaAn}.\eop

The convexity requirement in part (a) amounts in Examples \ref{ex:goal} and \ref{ex:risk} to $G$ having convex component functions because $y\geq 0$ for $y\in Y$. In these examples, $Y$ is also compact and thus (a) applies. 

If the condition in Theorem \ref{thm:duality}(b) holds with $\inf_{x\in \reals^n} \phi(x)$ finite, then $\sup_{y\in \reals^m} \psi(y)$ is also finite and there exist $x\in \reals^n$ and $y\in \reals^m$ that attain the minimum value of $\phi$ and maximum value of $\psi$, respectively. Likewise, if $\sup_{y\in \reals^m} \psi(y)$ is finite, then $\inf_{x\in \reals^n} \phi(x)$ is also finite and the minimum value and the maximum value are again attained; see \cite[Theorem 11.42]{VaAn}. We observe that this is not the case under Theorem \ref{thm:duality}(a). As a counterexample, suppose that $X = [-1,1]$, $G(x) = (x, x^2)$, and $h(z) = \sup_{y\in \{1\} \times [0,\infty)} \langle y,z\rangle = z_1 + \iota_{(-\infty,0]}(z_2)$, which then fits Theorem \ref{thm:duality}(a). Elementary calculations lead to 
\[
\inf_{x\in \reals} \phi(x) = 0 = \sup_{y\in \reals^2} \psi(y)  ~~\mbox{ because } ~~\psi(y) = \begin{cases}
-\infty & \mbox{ if } y \not\in \{1\} \times [0,\infty)\\
-1 + y_2 & \mbox{ if } y_1 = 1, y_2 \in [0,1/2)\\
-\tfrac{1}{4y_2} & \mbox{ if } y_1 = 1, y_2 \in [1/2,\infty).
\end{cases}   
\]
However, the maximum value of $\psi$ is not attained.

Without convexity or linearity in $G$, we effectively need to move beyond Theorem \ref{thm:duality} and the Rockafellian defined in \eqref{eqn:Rock} to close the likely gap between the maximum value of $\psi$ and the minimum value of $\phi$. The following general condition provides a geometric interpretation of what is at stake. 

\begin{theorem}{\rm (exactness).}\label{thm:exactness}  
Suppose that $f:\reals^q\times \reals^n\to \Reals$ is a Rockafellian for $g:\reals^n\to\Reals$. If there exists $\bar y\in \reals^q$ such that    
\begin{equation}\label{eqn:defexactness}
\inf_{x\in \reals^n} f(u,x) \geq \inf_{x\in \reals^n} g(x) + \langle \bar y, u\rangle ~~~\forall u\in\reals^q, 
\end{equation}
then the resulting dual function given by $\psi(y) = \inf_{u,x} f(u,x) - \langle y,u\rangle$ satisfies 
\[
\inf_{x\in \reals^n} g(x) = \psi(\bar y) = \sup_{y\in \reals^q} \psi(y).    
\]
\end{theorem}
\state Proof. Trivially, $\psi(y) \leq \inf g$ for all $y\in \reals^q$ because $\inf_{u \in \reals^q} f(u,x) - \langle y, u\rangle \leq f(0,x) = g(x)$ regardless of $x$ by virtue of $f$ being a Rockafellian for $g$. The assumption \eqref{eqn:defexactness} amounts to having $\psi(\bar y) = \inf_{u,x} f(u,x) - \langle \bar y,u\rangle \geq \inf g$.\eop

The theorem states that when $\inf g$ is finite, then the ability to construct a supporting hyperplane at the point $(0, \inf g)$ to the epigraph of the min-value function $u\mapsto \inf_{x\in \reals^n} f(u,x)$ is the key property. It can be satisfied by constructing Rockafellians that penalize nonzero $u$-values ``sufficiently.'' 

Returning to the setting of Section \ref{sec:intro}, the insight from Theorem \ref{thm:exactness} brings us to a second Rockafellian for $\phi$ in the actual problem \eqref{eqn:actualproblem2}: 
\begin{equation}\label{eqn:Rock2}
f_\theta(u,x) = \iota_X(x) + h\big(G(x) + u\big) + \tfrac{1}{2}\theta\|u\|_2^2,
\end{equation}
where $\theta \in (0,\infty)$ is a parameter in a quadratic term not present in the first Rockafellian $f$; see \eqref{eqn:Rock}. Certainly $f_\theta(u,x) \geq f(u,x)$ for all $u,x$ and the difference between the Rockafellians might be sufficiently large to ensure exactness in the sense of Theorem \ref{thm:exactness}. Regardless, the Lagrangian resulting from $f_\theta$ takes the form (see, e.g., \cite[Example 6.7]{primer})
\begin{equation}\label{eqn:augLagr}
l_\theta(x,y) = \iota_X(x) - \inf_{w\in Y} \Big\{ \tfrac{1}{2} \langle w, Qw\rangle - \big\langle G(x), w\big\rangle + \tfrac{1}{2\theta} \big\|w - y\big\|_2^2\Big\},
\end{equation}
which in turn defines a new dual problem:
\[
\nnmax_{y\in \reals^m} \, \psi_\theta(y) = \inf_{x\in\reals^n} l_\theta(x,y).
\]
In view of its construction from the first Rockafellian $f$ by adding a term $\tfrac{1}{2}\theta\|u\|_2^2$, the second Lagrangian $l_\theta$ is commonly called an {\em augmented Lagrangian.}

The verification of Theorem \ref{thm:exactness} for $f_\theta$ is supported by the following fact. Suppose that $X$ is compact and $\bar y\in \reals^m$. If there are $\bar\theta\in [0,\infty)$ and a neighborhood $U$ of $0\in\reals^m$ such that 
\begin{equation}\label{eqn:augLagrCond}
\inf_{x\in \reals^n} f(u,x) \geq \inf_{x\in \reals^n} \phi(x) + \langle \bar y, u\rangle - \tfrac{1}{2}\bar\theta \|u\|_2^2~~~\forall u\in U,
\end{equation}
then, for sufficiently large $\theta$,  
\[
\inf_{x\in \reals^n} f_{\theta}(u,x) \geq \inf_{x\in \reals^n} \phi(x) + \langle \bar y, u\rangle ~~~\forall u\in\reals^m. 
\]
The key benefit from this fact, which follows directly because $u\mapsto\inf_{x\in \reals^n} f(u,x)$ is bounded from below by an affine function, is that we only need to consider a neighborhood $U$ in \eqref{eqn:augLagrCond} as compared to every $u$ in Theorem \ref{thm:exactness}.   

\begin{example}{\rm (nonlinear optimization (cont.)).}\label{ex:nonlinear2} In the setting of Example \ref{ex:nonlinear} but with no inequality constraints, 
the Rockafellian in \eqref{eqn:Rock} produces the Lagrangian
\[
l(x,y) = g_0(x) + \sum_{i=1}^m y_i g_i(x) ~~\mbox{ when }~y = (1, y_1, \dots, y_m).
\]
The Rockafellian in \eqref{eqn:Rock2} defines the Lagrangian
\[
l_\theta(x,y) = g_0(x) + \sum_{i=1}^m y_i g_i(x) + \tfrac{1}{2}\theta \sum_{i=1}^m \big(g_i(x)\big)^2 ~~\mbox{ when }~y = (1, y_1, \dots, y_m).
\]
The simple instance with $g_0(x) = -x^2$ and $g_1(x) = x$ reveals the difference between these Lagrangians. The first Rockafellian produces $\inf_x f(u,x) = u_0 -u_1^2$ so the requirement of Theorem \ref{thm:exactness} cannot hold. In fact, the dual function has $\psi(y) = -\infty$ for all $y\in\reals^2$ while $\inf \phi = 0$. The second Rockafellian has  $\inf_x f_\theta(u,x) = u_0-u_1^2 + (1/2)\theta (u_0^2+u_1^2)$. Thus, for $\theta\in [2,\infty)$ one can employ $\bar y = (1,0)$ to satisfy the requirement of Theorem \ref{thm:exactness} and $\psi_\theta(\bar y) = \inf \phi$ for such $\theta$.
\end{example}

An additional advantage of the augmented Lagrangian in \eqref{eqn:augLagr} is that $l_\theta(x, \cdot)$ is smooth for $x\in X$, while $l(x,\cdot)$ in \eqref{eqn:lagr} may not have that property. We can see this by examining the inf-term in \eqref{eqn:augLagr}, which amounts to a Moreau envelope of the convex function 
\[
y\mapsto \iota_Y(y) + \tfrac{1}{2} \langle y, Qy\rangle - \big\langle G(x), y\big\rangle. 
\]
Its smoothness follows by \cite[Theorem 2.26]{VaAn}.

Dual problems and the corresponding Rockafellian and Lagrangian relaxations often turn out to be better behaved than the original problem. They absorb inaccuracies in problem data more easily resulting in convergence properties of the kind seen in Section \ref{sec:approx} but under milder conditions \cite{DerideRoyset.24}. Regardless of the Rockafellian, the resulting dual function $\psi$ is concave (i.e., $-\psi$ is convex), which certainly is beneficial. 

A Lagrangian sets up a game between an $x$-player that seeks to minimize $l(\cdot,y)$ and a $y$-player that aims to maximize $l(x,\cdot)$. This is fertile ground for algorithms that alternate between these two players in an effort to identify a {\em saddle point} and thus simultaneously solve the original problem and its dual. The optimality condition in Theorem \ref{thm:OptimComposite} can be viewed from this angle. Adopting the Rockafellian in \eqref{eqn:Rock} and the resulting Lagrangian from \eqref{eqn:lagr}, the optimality condition is equivalently stated as 
\[
0 \in \partial_x l(x,y)~~~\mbox{ and }~~~ 0\in \partial_y(-l)(x,y)
\]   
because the first inclusion simplifies to $0 \in N_X(x) + \nabla G(x)^\top y$ and the second inclusion amounts to $0 \in -G(x) + Qy + N_Y(y)$. When $G$ is affine, then the actual problem is convex and these inclusions are equivalent to a saddle point condition for $l$; see, e.g., \cite[Proposition 5.36]{primer} and \cite[Section 5.E]{primer}. 

A shift to the Rockafellian in \eqref{eqn:Rock2} can be viewed as a means to elicit a saddle point in a local sense and thus bypass the need for convexity globally. This results in {\em augmented Lagrangian methods} that essentially alternate between minimizing $l_{\theta^\nu}(\cdot,y^\nu)$, at least locally near a current point $x^\nu$, for a multiplier vector $y^\nu$ to produce $x^{\nu+1}$ and updating 
\[
y^{\nu+1} = y^\nu + \lambda^\nu \nabla_y l_{\theta^\nu}(x^{\nu+1}, y^\nu) = y^\nu + \frac{\lambda^\nu}{\theta^\nu}(\hat w - y), 
\]
where $\lambda^\nu$ is a positive step size coordinated with the penalty parameter $\theta^\nu$. Both quantities may change across the iterations. The $y$-update leverages the fact that $l_\theta(x,\cdot)$ is smooth as noted earlier, and this produces the stated formula with $\hat w$ being the minimizer in \eqref{eqn:augLagr} under $x = x^{\nu+1}$. Recent developments in this direction appear in \cite{Rockafellar.23b,HangSarabi.23b}, which also include historical remarks about such methods. The analysis of augmented Lagrangian methods requires second-order variational analysis, which we briefly summarize next.

\section{Second-Order Theory}\label{sec:secondorder}

We know from classical analysis of twice smooth functions that second-order differentiability is central when refining a first-order optimality condition: $\nabla f(\bar x) = 0$ together with a positive definite Hessian $\nabla^2 f(\bar x)$ ensure that $\bar x$ is a ``stable'' local minimizer of $f$ in a specific sense. Concretely, suppose that  
\begin{equation}\label{eqn:quadex}
f(x) = \tfrac{1}{2}\langle x, Bx\rangle + \langle c, x\rangle, 
\end{equation}
where $B$ is a symmetric positive definite matrix. The first-order condition $\nabla f(x) = 0$ has the unique solution $\bar x = -B^{-1}c$, with the solution remaining unique and changing {\em proportionally} when the right-hand side is modified from 0 to some vector $y$. In detail, if $s(y) = \{x \, | \, \nabla f(x) = y\}$, then 
\begin{equation}\label{eqn:sLip}  
\big\|s(y) - s(\bar y)\big\|_2 \leq \|B^{-1}\| \|y-\bar y\|_2
\end{equation}
for any matrix norm $\|\cdot\|$ consistent with the Euclidean norm. Consequently, the solution of $\nabla f(x) = y$ is unique and Lipschitz continuous as a function of $y$. A nonzero residual $y$ for a candidate solution $x'$ of the equation $\nabla f(x) = 0$, which is unavoidable in numerical computations, would therefore be of minor concern because the distance from $x'$ to the actual solution $s(0)$ is at most $\|B^{-1}\| \|y\|_2$. The equation $\nabla f(x) = 0$ is well posed. The same holds for the problem of minimizing $f$ due to the Fermat rule, which is both sufficient and necessary in this convex case. 

The situation is different if $f(x) = 0$ when $x<0$ and $f(x) = x^4$ otherwise. Then, $x = (y/4)^{1/3}$ solves $\nabla f(x) = y$ if $y\geq 0$ but there are no solutions if $y<0$. The problem of minimizing $f$ is then ill-posed in some sense. A difference between the two cases is revealed by the Hessian matrices. In the first case, $\nabla^2 f(\bar x)$ is positive definite. In the second case, the Hessian matrices are merely positive semidefinite at the minimizers. These elementary examples recall how the properties of a Hessian matrix can help characterize the stability of local minimizers of twice smooth functions. Variational analysis furnishes far reaching extensions for nonsmooth functions.

With Hessian matrices playing such a central role in stability analysis as well as in algorithms such as Newton's method, a goal is to define an analogous concept in the absence of smoothness and we follow the pioneering work in \cite{Mordukhovich.80}. For a set-valued mapping $S:\reals^n\tto \reals^m$, its {\em coderivative} at $(\bar x, \bar y) \in \gph S$ is the set-valued mapping $D^*S(\bar x, \bar y):\reals^m \tto \reals^n$ given by   
\[
D^*S(\bar x, \bar y)(v) = \big\{u\in \reals^n~\big|~(u,-v) \in N_{\gph S}(\bar x, \bar y) \big\} ~~~\mbox{ for } v\in \reals^m. 
\]

\begin{figure}
\centering
\includegraphics[width=0.425\textwidth]{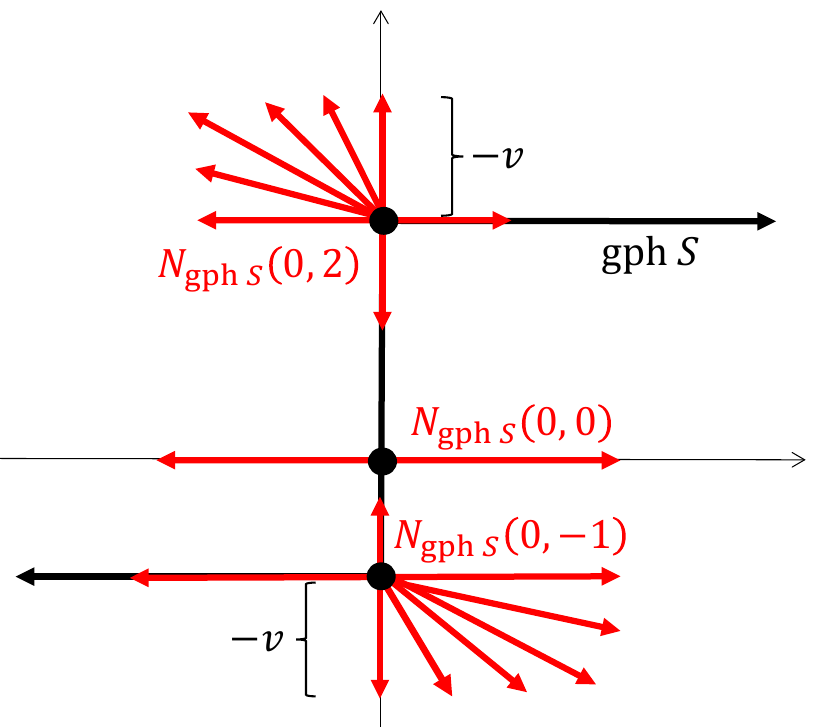}
\caption{Normal cones of $\gph S$ at the points $(0,-1)$, $(0,0)$, and $(0,2)$ define coderivatives.}\label{fig:coderiv}
\end{figure}

Figure \ref{fig:coderiv} illustrates the definition using the set-valued mapping $S:\reals\tto\reals$ with $S(x) = \{-1\}$ if $x<0$, $S(0) = [-1,2]$, and $S(x) = \{2\}$ if $x>0$. Near $(\bar x, \bar y) = (0,0)$, $\gph S$ is a vertical line and thus
\[
D^*S(0,0)(v) = \begin{cases}
  \reals & \mbox{ if } v=0\\
  \emptyset & \mbox{ otherwise. }
\end{cases} 
\]
At $(\bar x, \bar y) = (0,-1)$, there are normal vectors fanning out across the fourth quadrant and also in the other axis directions. In fact, $N_{\gph S}(0,-1)$ is a nonconvex cone producing the coderivative:
\[
D^*S(0,-1)(v) = \begin{cases}
  \reals     & \mbox{ if } v = 0\\
  [0,\infty) & \mbox{ if } v\in (0,\infty)\\
  \{0\} & \mbox{ otherwise. }
\end{cases} 
\]
At $(\bar x, \bar y) = (0,2)$, the normal cone is also nonconvex and coincides with the second quadrant augmented with the axis directions. Consequently, we obtain 
\[
D^*S(0,2)(v) = \begin{cases}
  \reals & \mbox{ if } v = 0\\
  (-\infty, 0] & \mbox{ if } v\in (-\infty, 0)\\
  \{0\} & \mbox{ otherwise. }
\end{cases} 
\]

We are especially interested in the case when $S(x) = \partial f(x)$ for some function $f:\reals^n\to \Reals$. This leads to the key definition from \cite{Mordukhovich.92}: For $f:\reals^n\to \Reals$ and a point $\bar x$ where $f(\bar x)$ is finite, the {\em second-order subdifferential} of $f$ at $\bar x$ relative to $\bar y\in \partial f(\bar x)$ is defined as 
\[
\partial^2 f(\bar x, \bar y)(v) = D^*(\partial f)(\bar x,\bar y)(v)~~~~\mbox{ for } v\in \reals^n.
\] 

The elementary example $f(x) = x^2$ produces $\partial^2 f(\bar x, \bar y)(v) = 2v$ at a point $(\bar x, f(\bar x))$ because normal vectors at any point of the graph of $x\mapsto 2x$ are of the form $\lambda(2,-1)$ with $\lambda \in \reals$. More generally, if $f:\reals^n\to \Reals$ is twice smooth in a neighborhood of $\bar x$, then 
\begin{equation}\label{eqn:partial2C2}
\partial^2 f\big(\bar x, \nabla f(\bar x)\big)(v) = \nabla^2 f(\bar x) v ~~~~\mbox{ for } v\in \reals^n.
\end{equation}
If $f$ is merely smooth with locally Lipschitz continuous gradients in a neighborhood of $\bar x$, then 
\[
\partial^2 f\big(\bar x, \nabla f(\bar x)\big)(v) = \partial g_v(\bar x)~~~~\mbox{ for } v\in \reals^n, ~~\mbox{ where } g_v(x) = \big\langle v, \nabla f(x)\big\rangle.
\]
An example of the latter situation is furnished by $f(x) = x^2$ if $x >0$ and $f(x) = 0$ otherwise. Here, $g_v(x) = 2vx$ if $x > 0$ and $g_v(x) = 0$ otherwise. Thus,  
\[
\partial^2 f\big(\bar x, \nabla f(\bar x)\big)(v) = \begin{cases}
 \{2v\} & \mbox{ if } x > 0\\ 
 [0, 2v]  & \mbox{ if } v\geq 0, x = 0\\ 
 \{0, 2v\}  & \mbox{ if } v< 0, x = 0\\ 
 \{0\}  & \mbox{ if } x < 0. 
\end{cases}
\]

With the second-order subdifferential for arbitrary functions now well defined, we return to the question of solution stability, which first needs to be formalized. We follow \cite{PoliquinRockafellar.98} and say that a point $\bar x$ is a {\em tilt-stable local minimizer} of $f:\reals^n\to \Reals$ when $f(\bar x)$ is finite and there exists $\delta >0$ such that the set-valued mapping $M:\reals^n\tto\reals^n$ given by 
\begin{equation}\label{eqn:Mmap}
M(y) = \nargmin_x \big\{ f(x) - f(\bar x) - \langle y,x\rangle~\big|~ \|x - \bar x\|_2 \leq \delta \big\}
\end{equation}
is single-valued and Lipschitz on some neighborhood of $y = 0$, and has $M(0) = \bar x$. Here, the Lipschitz property amounts to having a finite $\kappa$ such that $\|M(y'') - M(y')\|_2 \leq \kappa \|y'' - y'\|_2$ for all $y''$ and $y'$ sufficiently close to the origin. 

For a tilt-stable local minimizer $\bar x$ of $f$, it is immediately clear from the Fermat rule that $M(y) \subset \{x \, | \, y \in \partial f(x)\}$ when $y$ is sufficiently close to the origin. The gradual change of a solution to an optimality condition under different tolerances---the motivation behind the initial discussion around \eqref{eqn:quadex}---is therefore closely related to tilt-stability with the Lipschitz property emerging concretely in \eqref{eqn:sLip}. We clarify this point using two facts.

First, in the case of a twice smooth function, a tilt-stable local minimizer at $\bar x$ is equivalent to having a positive definite Hessian matrix at $\bar x$ and thus the discussion around \eqref{eqn:quadex} could as well have used the terminology of tilt-stability. Formally, \cite[Proposition 1.2]{PoliquinRockafellar.98} states:   

\begin{proposition}{\rm (tilt-stability under twice smoothness).}\label{thm:tiltsmooth} 
For a twice smooth function $f:\reals^n\to \reals$ and a point $\bar x$ with $\nabla f(\bar x) = 0$ one has: 
\[
\bar x \mbox{ is a tilt-stable local minimizer of $f$}  ~~~\Longleftrightarrow~~~  \nabla^2 f(\bar x) \mbox{ is positive definite}.
\]
\end{proposition}

Second, in the case of a general function, we need two additional concepts. A function $f:\reals^n\to \Reals$ is {\em prox-regular} at $\bar x$ for $\bar y$ if $\epi f$ is closed relative to a neighborhood of $(\bar x, f(\bar x))$, $\bar y\in \partial f(\bar x)$, and there exist $\theta,\epsilon>0$ such that 
\begin{align*}
  & f(x') > f(x) + \langle y, x' - x\rangle - \tfrac{1}{2}\theta \|x' - x\|_2^2 ~\mbox{ for } x'\neq x \mbox{ when }\\
  & \|x' - \bar x\|_2 < \epsilon, ~\|x - \bar x\|_2 < \epsilon, ~\big|f(x) - f(\bar x)\big| < \epsilon, ~\|y-\bar y\|_2 < \epsilon, ~y\in \partial f(x).  
\end{align*}
The function is {\em subdifferentially continuous} at $\bar x$ for $\bar y$ if $\bar y\in \partial f(\bar x)$ and, whenever $(x^\nu,y^\nu) \to (\bar x, \bar y)$ with $y^\nu\in \partial f(x^\nu)$, one has $f(x^\nu)\to f(\bar x)$.

Every proper, lsc, and convex function is prox-regular and subdifferentially continuous at any point of its domain for any subgradient at that point \cite[Example 13.30]{VaAn}. In the setting of Section \ref{sec:intro} with $G$ being twice smooth, the function given by $f(x) = h(G(x))$ is prox-regular and subdifferentially continuous at $\bar x$ for $\bar y\in \partial f(\bar x)$ if the qualification \eqref{eqn:compositeQual} holds \cite[Proposition 13.32]{VaAn}. However, $x\mapsto \min\{0,x\}$ is subdifferentially continuous but not prox-regular at zero; see \cite[Section 13.F]{VaAn} for further details.  

We are then ready to state the second fact about tilt-stability \cite[Theorem 1.3]{PoliquinRockafellar.98}: 

\begin{theorem}{\rm (tilt-stability).}\label{thm:tilt} 
For a function $f:\reals^n\to \Reals$ and a point $\bar x$ with $0\in \partial f(\bar x)$, suppose that $f$ is prox-regular and subdifferentially continuous at $\bar x$ for $\bar y = 0$. Then, one has  
\[
\bar x \mbox{ is a tilt-stable local minimizer of } f ~~\Longleftrightarrow ~~\langle u, v\rangle > 0 \mbox{ whenever } v \neq 0, u\in \partial^2 f(\bar x, 0)(v).
\]
These properties imply the existence of $\delta>0$ such that $M$ from \eqref{eqn:Mmap} has $M(y) = \{x~|~y\in \partial f(x)\}$ for all $y$ in a neighborhood of 0.
\end{theorem}

If $f$ is twice smooth, then it follows via \eqref{eqn:partial2C2} that the condition involving $\langle u, v\rangle >0$ in the theorem simplifies to 
\[
\big\langle \nabla^2 f(\bar x)v, v\big\rangle > 0 ~\mbox{ whenever } v\neq 0,
\]
which, as expected from Theorem \ref{thm:tiltsmooth}, is equivalent to $\nabla^2 f(\bar x)$ being positive definite. 

The example $f(x) = \max\{-x,2x\}$ has $\partial f(x) = S(x)$ from Figure \ref{fig:coderiv}. We can use Theorem \ref{thm:tilt} to check whether $\bar x = 0$ is a tilt-stable local minimizer of $f$. Since $f$ is proper, lsc, and convex, it is prox-regular and subdifferentially continuous at $\bar x$ for $\bar y = 0$. Figure \ref{fig:coderiv} shows that $\partial^2 f(0,0)(v) = \emptyset$ for $v\neq 0$. Thus, $\bar x = 0$ is a tilt-stable local minimizer of $f$. 

The situation is different for $f(x) = \max\{0,x\}$ at $\bar x = 0$. While certainly a (local) minimizer, $\bar x$ is not tilt-stable because $\partial^2 f(0,0)(v) = [0,\infty)$ when $v>0$ and   $\partial^2 f(0,0)(v) = \{0\}$ when $v< 0$. Thus, in either case, one can select $u = 0$ in Theorem \ref{thm:tilt}.

Returning to the setting of Section \ref{sec:intro}, we recall a characterization of tilt-stability \cite[Theorem 5.4]{MordukhovichRockafellar.12}:

\begin{theorem}{\rm (tilt-stability for actual problem).}\label{thm:tilt2} For the function $f:\reals^n\to \Reals$ given by 
\[
f(x) = h\big(G(x)\big),
\]
with $h$ and $G$ as defined in Section \ref{sec:intro}, and a point $\bar x\in\dom f$, suppose that $G$ is twice smooth, $Q$ defining $h$ is either positive definite or the zero matrix, the qualification \eqref{eqn:compositeQual} holds at $\bar x$, and $\bar y$ is a unique vector satisfying 
\[
\bar y \in \nargmin_{y\in Y} \tfrac{1}{2}\langle y, Qy\rangle - \big\langle y, G(\bar x)\big\rangle ~~\mbox{ and } ~~ \nabla G(\bar x)^\top \bar y  = 0. 
\]
Moreover, suppose that the following second-order qualification holds: 
\[
u \in \partial^2 h\big(G(\bar x),\bar y\big)(0) ~~\mbox{ and }~~ \nabla G(\bar x)^\top u = 0 ~~~\Longrightarrow ~~~  u = 0. 
\]
Then, with the notation $g(x) = \langle \bar y, G(x)\rangle$, one has 
\begin{align*}
\mbox{$\bar x$ is a tilt-stable local minimizer of $f$} ~~\Longleftrightarrow~~~ & \langle u, v\rangle > 0 ~\mbox{ whenever }~ v \neq 0,\\
& u \in \nabla^2g(\bar x) v + \nabla G(\bar x)^\top \partial^2 h\big(G(\bar x), \bar y\big)\big(\nabla G(\bar x)v\big).
\end{align*}
\end{theorem}

The proof of \cite[Theorem 5.4]{MordukhovichRockafellar.12} leverages the extensive calculus for second-order subdifferentials. It allows us to focus on the second-order subdifferential of $h$, which presumably is somewhat simpler than that of $f$. Specifically, \cite[Lemma 4.4]{MordukhovichRockafellar.12} and its proof reveal that
\[
u \in \partial^2 h(z,y)(v) ~\Longleftrightarrow~ Qu - v \in \partial^2 \iota_Y(y, z-Qy)(-u).
\]
Since $Y$ is a polyhedral set, the second-order subdifferential of $\iota_Y$ is precisely characterized in terms of its faces; see \cite[Lemma 4.4 ]{MordukhovichRockafellar.12} for details. As an example, $Y = [0, \infty)$ results in  $\gph N_{Y} = (\{0\} \times (-\infty,0]) \cup ([0,\infty) \times \{0\})$. Thus, with $C = \gph N_{Y}$, we obtain    
\[
\partial^2 \iota_{Y}(y,w)(v) = D^*(N_Y)(y,w)(v) = \big\{ u \, \big| \, (u,-v) \in N_C(y,w) \big\} = \begin{cases}
  \{0\} & \mbox{ if } y > 0, w = 0, v\in\reals\\
  (-\infty, 0] & \mbox{ if } y=0, w = 0, v<0\\
  \reals & \mbox{ if } y=0, w \leq 0, v=0\\
  \{0\} & \mbox{ if } y=0, w = 0, v>0\\
  \emptyset & \mbox{ otherwise.}   
\end{cases}
\]
Additional formulas for second-order subdifferentials of indicator functions of polyhedral sets as well as of $h$ when $Q$ is the zero matrix appear in \cite{MordukhovichSarabi.16, MordukhovichSarabi.16b}.

Second-order subdifferentials stand out with their versatile calculus as demonstrated, for example, in \cite{MordukhovichRockafellarSarabi.13,MordukhovichSarabi.16,MordukhovichSarabi.17,MohammadiMordukhovichSarabi.22}. In particular, they enter as means to characterize Lipschitz-like properties of set-valued mappings representing optimality conditions in composite optimization \cite{DoMordukhovichSarabi.21}. However, second-order theory is much richer as can be seen from the monographs \cite[Chapter 13]{VaAn} and, especially, \cite{Mordukhovich.24}. Developments include those leveraging twice epi-differentiability \cite{MohammadiSarabi.20,HangSarabi.24}. Second-order optimality conditions emerge as well; see \cite{MohammadiSarabi.20,MohammadiMordukhovichSarabi.22} for recent advances extending the classical work in \cite{Rockafellar.89,BurkePoliquin.92}. Current trends center on {\em variational convexity} as the key concept \cite{Rockafellar.23,KhanhKhoaMordukhovichPhat.24}. 

Perturbations beyond ``tilting'' are of major concerns as seen from the discussion in Section \ref{sec:approx}. They lead to the question of {\em full stability} of parametric optimization and generalized equations \cite{LevyPoliquinRockafellar.00}, which often can be addressed using second-order subdifferentials and related concepts \cite{MordukhovichRockafellarSarabi.13,MordukhovichNghiaRockafellar.15,MordukhovichSarabi.17,BenkoRockafellar.24}. Second-order theory also gives rise to numerous algorithms. Algorithmic development leveraging advanced theories from variational analysis is currently an active area of research. Some recent efforts include extensions of sequential quadratic programming \cite{Sarabi.22} and Newton-type methods \cite{MordukhovichSarabi.20,BurkeEngle.20,AAMorPerezaros.24}.

\section*{Declarations}

This work is supported in part by the Office of Naval Research under grants N00014-24-1-2741, N00014-24-1-2318, and N00014-24-1-2492.\\ 

\noindent Data Availability: There is no data associated with this manuscript.\\

\noindent Competing interests: The author is an editor of the journal.

\bibliographystyle{plain}
\bibliography{refs}

\end{document}

%% file: mac.tex
\usepackage{theorem}
\usepackage{enumerate}
\usepackage{array}
\usepackage[reqno]{amsmath}
\usepackage{amssymb}
\usepackage{mathtools}
\usepackage{latexsym}
\usepackage{makeidx}
\usepackage{fancybox}
\input epsf.sty
\usepackage[dvips]{graphicx,color}
\usepackage{showlabels}
\usepackage{subcaption}
\usepackage[dvipsnames]{xcolor}
\usepackage[normalem]{ulem}

\theoremstyle{change}
\newtheorem{proclaim}{PROCLAIM}[section]
\newtheorem{theorem}[proclaim]{Theorem}

\newtheorem{proposition}[proclaim]{Proposition}

\newtheorem{example}[proclaim]{Example}

\numberwithin{equation}{section}

\outer\def\proclaim #1. #2\par{\medbreak \noindent{\bf#1.\enspace}{\sl#2}\par
  \ifdim\lastskip<\medskipamount
  \removelastskip\penalty55\medskip\fi}
\def\state #1. { \noindent{\bf#1.\enspace}}
\def\algo #1. { \noindent{\bf#1.\enspace}}

\DeclareMathOperator{\dist}{dist}

\DeclareMathOperator{\dom}{dom}
\DeclareMathOperator{\gph}{gph}
\DeclareMathOperator{\epi}{epi}

\newcommand{\comp}{\,{\raise 1pt \hbox{$\scriptstyle\circ$}}\,}

\newcommand{\reals}{\mathbb{R}}
\newcommand{\Reals}{\overline{\mathbb{R}}}
\newcommand{\natnums}{{{\rm l} \kern -.13em {\rm N} }}
\newcommand{\nats}{\mathbb{N}}
\newcommand{\snats}{{I\kern -.29em N}}
\newcommand{\rats}{{Q\kern -.64em \raise 1pt \hbox{$\scriptstyle |$}\;\,}}
\newcommand{\srats}
	{{Q\kern -.56em \raise 1.2pt \hbox{$\scriptscriptstyle /$}\,}}
\newcommand{\ints}{Z\kern -.46em Z}

\newcommand{\pluss}{\hskip1pt \raise1pt\vbox{\hrule width6pt \vskip1pt \hrule
                    width6pt} \kern-4pt{\lower1pt\hbox{\vrule height6pt
		    \kern1pt\vrule height6pt}}\hskip5pt}
\newcommand{\eop}
	{\hfill{$\vcenter{\hrule height1pt \hbox{\vrule width1pt height5pt
   	 \kern5pt \vrule width1pt} \hrule height1pt}$} \medskip}

\newcommand{\setd}{{ d \kern -.15em l}}
\newcommand{\hatsetd}{ d \hat{\kern -.15em l }}

\renewcommand{\epsilon}{\varepsilon}
\renewcommand{\phi}{\varphi}

\newcommand{\tto}{\;{\lower 1pt \hbox{$\rightarrow$}}\kern -12pt
           \hbox{\raise 2.5pt \hbox{$\rightarrow$}}\;}
\newcommand{\overto}[1]{\,{\raise 0pt\hbox{$\rightarrow$}}\kern -9pt
     \hbox{\lower 3pt \hbox{$\scriptscriptstyle#1$}}\hskip6pt}
\newcommand{\underto}[1]{\,{\lower 1pt\hbox{$\rightarrow$}}\kern -9pt
     \hbox{\raise 4pt \hbox{$\,\scriptscriptstyle#1$}}\hskip7pt}
\newcommand{\bigoverto}[1]{{\raise 0pt\hbox{$\,\longrightarrow$}}\kern -16pt
     \hbox{\lower 3pt \hbox{$\scriptscriptstyle#1$}}\hskip4pt}
\newcommand{\bigunderto}[1]{\,{\lower 1pt\hbox{$\longrightarrow$}}\kern -16pt
     \hbox{\raise 4pt \hbox{$\,\scriptscriptstyle#1$}}\hskip6pt}
\newcommand{\bigbigto}[2]{\,{\raise 0pt\hbox{$\,\longrightarrow$}}\kern -16pt
     \hbox{\lower 3pt \hbox{$\scriptscriptstyle#2$}}\kern -10pt
     \hbox{\raise 4pt \hbox{$\,\scriptscriptstyle#1$}}\hskip7pt}
\newcommand{\downto}{{\raise 1pt \hbox{$\scriptscriptstyle \,\searrow\,$}}}
\newcommand{\upto}{{\raise 1pt \hbox{$\scriptscriptstyle \,\nearrow\,$}}}

\newcommand{\notimply}
	{\quad\hbox{$\Longrightarrow \kern -14pt {/}$}\hskip6pt\quad}

\newcommand{\lto}{\,{\lower 1pt\hbox{$\rightarrow$}}\kern -10pt
     \hbox{\raise 4pt \hbox{$\, \scriptstyle l$}}\hskip7pt}
\newcommand{\eto}{\,{\lower 1pt\hbox{$\rightarrow$}}\kern -10pt
     \hbox{\raise 4pt \hbox{$\, \scriptstyle e$}}\hskip7pt}
\newcommand{\hto}{\,{\lower 1pt\hbox{$\rightarrow$}}\kern -11pt
     \hbox{\raise 4pt \hbox{$\, \scriptstyle h$}}\hskip7pt}
\newcommand{\pto}{\,{\lower 1pt\hbox{$\rightarrow$}}\kern -11pt
     \hbox{\raise 4.5pt \hbox{$\, \scriptstyle p$}}\hskip7pt}
\newcommand{\cto}{\,{\lower 1pt\hbox{$\rightarrow$}}\kern -11pt
     \hbox{\raise 4pt \hbox{$\, \scriptstyle c$}}\hskip7pt}
\newcommand{\gto}{\,{\lower 1pt\hbox{$\rightarrow$}}\kern -11pt
     \hbox{\raise 4.5pt \hbox{$\, \scriptstyle g$}}\hskip7pt}
\newcommand{\sto}{\,{\lower 1pt\hbox{$\rightarrow$}}\kern -10pt
     \hbox{\raise 4pt \hbox{$\, \scriptstyle s$}}\hskip7pt}
\newcommand{\awto}{\,{\lower 1pt\hbox{$\rightarrow$}}\kern -15pt
     \hbox{\raise 4pt \hbox{$\, \scriptstyle aw$}}\hskip7pt}
\def\Nto{\,{\raise 1pt\hbox{$\rightarrow$}}\kern -13pt
     \hbox{\lower 3pt \hbox{$\, \scriptstyle N$}}\hskip7pt}
\def\Cto{\,{\raise 1pt\hbox{$\rightarrow$}}\kern -14pt
     \hbox{\lower 3pt \hbox{$\, \scriptstyle C$}}\hskip7pt}
\def\fto{\,{\raise 1pt\hbox{$\rightarrow$}}\kern -14pt
     \hbox{\lower 3pt \hbox{$\, \scriptstyle f$}}\hskip7pt}

\newcommand{\low}[1]{{\lower1pt \hbox{$\scriptstyle #1$}}}
\newcommand{\loww}[1]{{\lower2pt \hbox{$\scriptstyle #1$}}}
\newcommand{\high}[1]{{\raise1pt \hbox{$\scriptstyle #1$}}}

\newcommand{\nsum}{\mathop{\sum}\nolimits}

\newcommand{\nliminf}{\mathop{\rm liminf}\nolimits}

\renewcommand{\limsup}{\mathop{\rm limsup}}

\newcommand{\nsup}{\mathop{\rm sup}\nolimits}

\newcommand{\nnmin}{\mathop{\rm minimize}}
\newcommand{\nnmax}{\mathop{\rm maximize}}

\newcommand{\nargmin}{\mathop{\rm argmin}\nolimits}

\newcommand{\lwdy}[2]{\mathrel{\mathop
        {\raisebox{0.1ex}{\null$#1$}}{\hbox{\kern -1.0em
	{\raisebox{-0.8ex}{$\scriptstyle{\;\to #2}$}}}}}}
\newcommand{\lwwdy}[2]{\mathrel{\mathop
        {\raisebox{0.2ex}{\null$#1$}}{\hbox{\kern -1.0em
	{\raisebox{-1.1ex}{$\scriptstyle{\;\to #2}$}}}}}}
\newcommand{\slwdy}[2]{\scriptsize{{\mathrel{\mathop
        {\raisebox{0.1ex}{\null$#1$}}{\hbox{\kern -1.0em
	{\raisebox{-0.8ex}{$\scriptstyle{\;\to #2}$}}}}}}}}
\newcommand{\slwwdy}[2]{\scriptsize{{\mathrel{\mathop
        {\raisebox{0.2ex}{\null$#1$}}{\hbox{\kern -1.0em
	{\raisebox{-1.1ex}{$\scriptstyle{\;\to #2}$}}}}}}}}

\definecolor{lightgray}{gray}{0.75}
\definecolor{myred}{rgb}{0.55,0,0}
\definecolor{myblue}{rgb}{0,0,0.5} 
\definecolor{mygreen}{rgb}{0,0.5,0} 
\definecolor{purple}{rgb}{0.5,0,0.5} 
\definecolor{turq}{rgb}{0,0.805,0.816} 
\definecolor{maroon}{rgb}{0.51,0,0}
\definecolor{MAROON}{rgb}{0.51,0,0}
\definecolor{redor}{rgb}{0.78,0.078,0.078}
\definecolor{dgreen}{rgb}{0,0.3,0}

\newcommand{\bcdot}{\,{\raise .2ex \hbox{$\centerdot$}}\,}

\newcommand{\bbA}{\mathbb{A}}